\documentclass[a4paper,10pt]{amsart} 

\usepackage[ansinew]{inputenc} 
\usepackage[T1]{fontenc}
\usepackage[english]{babel} 

\usepackage{enumerate} 
\usepackage{booktabs}

\usepackage{amsmath} 
\usepackage{amsfonts} 
\usepackage{amssymb}
\usepackage{amsthm}
\usepackage{latexsym}
\usepackage{bbm}
\usepackage{a4wide}
\usepackage[pdftex]{graphicx}
  \DeclareGraphicsExtensions{.pdf,.png,.jpg,.jpeg,.mps}
  \usepackage{pgf}
  \usepackage{tikz}
\usepackage{floatflt}

\makeatletter

\@namedef{subjclassname@2010}{%

  \textup{2010} Mathematics Subject Classification}

\makeatother

\theoremstyle{plain} 
\newtheorem{theorem}{Theorem}[section] 
\newtheorem{lemma}[theorem]{Lemma} 
\newtheorem{prop}[theorem]{Proposition}

\newtheorem*{theorem*}{Theorem}

\theoremstyle{definition} 

\newtheorem*{definition*}{Definition}

\theoremstyle{remark}
\newtheorem{example}[theorem]{Example} 
\newtheorem*{remark}{Remark} 

\newcommand{\R}{\mathbbm{R}}

\newcommand{\Z}{\mathbbm{Z}}
\newcommand{\C}{\mathbbm{C}}
\newcommand{\E}{\mathbbm{E}}
\newcommand{\prob}{\mathbbm{P}}
\newcommand{\D}{\,\text{d}}

\title[Rademacher maximal function]{On the Rademacher maximal function}
\author[M. Kemppainen]{Mikko Kemppainen}
\address{Department of Mathematics and Statistics, University of Helsinki,
Gustaf Hällströmin katu 2b, FI-00014 Helsinki, Finland}

\email{mikko.k.kemppainen@helsinki.fi}

\begin{document}



\maketitle

\begin{abstract}

This paper studies a new maximal operator introduced by Hytönen, McIntosh and Portal in 2008
for functions taking values in a Banach space. The $L^p$-boundedness
of this operator depends on the range space; certain requirements on type and cotype are present for instance.
The original Euclidean definition of the maximal function is generalized to $\sigma$-finite measure spaces with filtrations
and the $L^p$-boundedness is shown not to depend on the underlying measure space or the filtration.
Martingale techniques are applied to prove that a weak type inequality is sufficient for $L^p$-boundedness and
also to provide a characterization by concave functions.

\end{abstract}

\tableofcontents

\section{Introduction}

The properties of the standard dyadic maximal function
\begin{equation*}
  Mf(\xi ) = \sup_{Q\ni\xi} | \langle f \rangle_Q | , \quad \xi\in\R^n ,
\end{equation*}
where $\langle f \rangle_Q$ denotes the average of a locally integrable function $f$ over a dyadic cube $Q$, are
well-known. More precisely, the (sublinear) operator $f\mapsto Mf$ is bounded in $L^p$ for all $p\in (1,\infty ]$
and satisfies for all $f\in L^1$ a certain weak type inequality (and also, is bounded from 
the dyadic Hardy space $H^1$ to $L^1$). These properties remain unchanged even if one studies functions
taking values in a Banach space and replaces absolute values by norms.

In their paper \cite{HMP}, Hytönen, McIntosh and Portal needed a new maximal function in order to prove
a vector-valued version of Carleson's embedding theorem. Instead of the supremum of (norms of) dyadic averages this
maximal function measures their R-bound (see Section 2 for the definition), 
which in general is not comparable to the supremum. More precisely, they
defined the Rademacher maximal function
\begin{equation*}
  M_Rf(\xi) = \mathcal{R} \Big( \langle f \rangle_Q : Q\ni\xi \Big) , \quad \xi\in\R^n ,
\end{equation*}
for functions $f$ taking values in a Banach space. They proved that the $L^p$-boundedness of $f\mapsto M_Rf$
is independent of $p$ in the sense that boundedness for one $p\in (1,\infty )$ implies boundedness for all 
$p$ in that range and that for many common range spaces including all UMD function lattices and spaces
with type 2, the operator $M_R$ is $L^p$-bounded. 
Nevertheless it turned out that the new maximal operator is not bounded for all choices of
range spaces, e.g. not for $\ell^1$.

The study of the Rademacher maximal operator continues here in a bit more general framework, which
was motivated by the need for vector-valued maximal function estimates in the context of
non-homogeneous spaces in \cite{HYTONENNONHOMTB}. We consider
it for operator-valued functions defined on 
$\sigma$-finite measure spaces, where averages are replaced by conditional expectations
with respect to filtrations. The boundedness of $M_R$ -- the RMF-property (of the range space) -- is shown 
not to depend on these new parameters;
instead, it is sufficient to check it for the filtration of dyadic intervals on $[0,1)$ (Theorem \ref{rmffilt}).
Here we follow
a reduction argument from Maurey \cite{MAUREYUMD}, originally tailored for the UMD-property. 
We also show that the RMF-property
requires non-trivial type and finite cotype of the Banach spaces involved
(Proposition \ref{nontrivialtype}). The Rademacher maximal function
is readily defined for martingales $X = (X_j)_{j=1}^{\infty}$ of operators by
\begin{equation*}
  X_R^* = \mathcal{R} \Big( X_j : j\in\Z_+ \Big) .
\end{equation*}
We will show using ideas from Burkholder \cite{BURKHOLDERUMD} that the RMF-property 
(requiring $L^p$-boundedness of $M_R$) is actually equivalent (Theorem \ref{weakrmf}) to the weak type
inequality (or the weak RMF-property)
\begin{equation*}
  \prob (X_R^* > \lambda ) \lesssim \frac{1}{\lambda} \| X \|_1 .
\end{equation*}
Finally, the RMF-property is characterized using concave functions
(Theorem \ref{concave}) in the spirit of Burkholder \cite{BURKHOLDEREXP}.

\section{Preliminaries}

All \emph{random variables} in Banach spaces (functions from a probability space to the Banach space) are
assumed to be $\prob$-strongly measurable, by which we mean that they are $\prob$-almost everywhere limits
of simple functions on the probability space whose measure we denote by $\prob$.
Their expectation, denoted by $\E$, is given by the Bochner integral. By an $L^p$-random
variable, for $1\leq p < \infty$, we mean a random variable $X$ (in a Banach space) whose $p$th moment 
$\E \| X \|^p$ is finite. 

Let $(\varepsilon_j)_{j=1}^{\infty}$ be a sequence of \emph{Rademacher variables}, i.e. a
sequence of independent random variables attaining values $-1$ and $1$ with an equal probability
$\prob (\varepsilon_j = -1) = \prob (\varepsilon_j = 1) = 1/2$. By the independence we have
$\E (\varepsilon_j \varepsilon_k) = (\E \varepsilon_j)(\E \varepsilon_k) = 0$, whenever $j\neq k$,
while (trivially) $\E (\varepsilon_j \varepsilon_k) = 1$, if $j=k$.
The equality of a randomized norm and a square sum of norms for vectors $x_1, \ldots , x_N$ in a Hilbert space
is thus established by the following calculation:
\begin{equation}
\label{randomsquare}
  \E \Big\| \sum_{j=1}^N \varepsilon_j x_j \Big\|^2 = \E \Big\langle \sum_{j=1}^N \varepsilon_j x_j ,
  \sum_{k=1}^N \varepsilon_k x_k \Big\rangle = \sum_{j,k=1}^N \E (\varepsilon_j \varepsilon_k)
  \langle x_j, x_k \rangle = \sum_{j=1}^N \| x_j \|^2 .
\end{equation}


The following standard result guarantees the comparability of different
randomized norms (see Kahane's book \cite{KAHANE} for a proof).

\begin{theorem} (The Khintchine-Kahane inequality)
\label{kkineq}
  For any $1 \leq p,q < \infty$, there exists a constant $K_{p,q}$ such that
  \begin{equation*}
    \Big( \E \Big\| \sum_{j=1}^N \varepsilon_j x_j \Big\|^p \Big)^{1/p}
    \leq K_{p,q} \Big( \E \Big\| \sum_{j=1}^N \varepsilon_j x_j \Big\|^q \Big)^{1/q},
  \end{equation*}
  whenever $x_1, \ldots , x_N$ are vectors in a Banach space.
\end{theorem}

The concepts of type and cotype of a Banach space
intend to measure how far the randomized norms are from $\ell^p$ sums of norms.

\begin{definition*}
  A Banach space $E$ is said to have 
  \begin{enumerate}
  \item \emph{type} $p$ for $1\leq p \leq 2$ if there exists a constant $C$ such that
  \begin{equation*}
    \Big( \E \Big\| \sum_{j=1}^N \varepsilon_j x_j \Big\|^2 \Big)^{1/2}
    \leq C \Big( \sum_{j=1}^N \| x_j \|^p \Big)^{1/p}
  \end{equation*}
  for any vectors $x_1, \ldots , x_N$ in $E$, regardless of $N$. 

  \item \emph{cotype} $q$ for $2\leq q \leq \infty$ if there exists a constant $C$ such that
  \begin{equation*}
    \Big( \sum_{j=1}^N \| x_j \|^q \Big)^{1/q}
    \leq C \Big( \E \Big\| \sum_{j=1}^N \varepsilon_j x_j \Big\|^2 \Big)^{1/2}
  \end{equation*} 
  for any vectors $x_1, \ldots , x_N$ in $E$, regardless of $N$. In the case $q=\infty$ the left hand side in the
  above inequality is replaced by $\max_{1\leq j \leq N} \| x_j \|$.

  \end{enumerate}
\end{definition*}

\begin{remark} A few observations are in order.
\begin{enumerate}
\item
As every Banach space has both type $1$ and cotype $\infty$ we say that
a Banach space has \emph{non-trivial type} (respectively \emph{finite cotype})
if it has type $p$ for some $p > 1$ (respectively cotype $q$ for some $q < \infty$).


\item
One can show that
$L^p$-spaces have type $\min \{ p,2 \}$ and cotype $\max \{ p,2 \}$ when $1\leq p < \infty$.
Sequence spaces $\ell^1$ and $\ell^{\infty}$ are, on the other hand, typical examples of spaces with only trivial type.

\item
Type and cotype of a Banach space $E$ and its dual $E^*$ are related in a natural way:
If $E$ has type $p$, then $E^*$ has cotype $p'$, where $p'$ is the Hölder conjugate of $p$.

\item
The equality (\ref{randomsquare}) of randomized norms and square sums of norms in Hilbert spaces means of course that
they have both type $2$ and cotype $2$.
A remarkable result of Kwapie\'n's (see the original paper \cite{KWAPIEN}, 
or the new proof by Pisier in \cite{PISIER}) 
is that a Banach space with both type $2$ and cotype $2$
is necessarily isomorphic to a Hilbert space.
\end{enumerate}
\end{remark}

The geometry of a Banach space can be studied by looking at 
its finite dimensional subspaces. 
We denote by $\ell_N^p$, where $p\in [1,\infty ]$ and $N\in\Z_+$, the $N$-dimensional 
subspace of $\ell^p$ consisting of sequences for which all but the
first $N$ terms are zero.
A Banach space $E$ is 
said to \emph{contain $\ell_N^p$'s $\lambda$-uniformly} for a $\lambda \geq 1$ if there exist for each $N\in\Z_+$ 
an $N$-dimensional subspace $E_N$ of $E$ and a bounded isomorphism $\Lambda_N : E_N \to \ell_N^p$ such that
$\| \Lambda_N \| \| \Lambda_N^{-1} \| \leq \lambda$.

The following theorem of Maurey and Pisier (see \cite{MP} for the original proof, or
\cite{DJT}, Theorems 13.3 and 14.1) relates this to the concept of type and cotype:

\begin{theorem}
\label{thm:contunif}
  Suppose that $E$ is a Banach space. Then
  \begin{enumerate}
    \item $E$ has a non-trivial type if and only if it does not contain $\ell_N^1$'s uniformly (i.e.
          $\lambda$-uniformly for some $\lambda \geq 1$).
    \item $E$ has finite cotype if and only if it does not contain $\ell_N^{\infty}$'s uniformly.
  \end{enumerate}
\end{theorem}

\begin{prop}
  If $E^*$ has non-trivial type, then $E$ has finite cotype.
  \begin{proof}
    Non-trivial type implies finite cotype for the dual and thus it follows from the assumption that
    $E^{**}$ has finite cotype. By Theorem \ref{thm:contunif}, $E^{**}$
    does not contain $\ell_N^{\infty}$'s uniformly and the same has to hold for its subspace $E$. This means
    that $E$ must have finite cotype.
  \end{proof}
\end{prop}

The proposition above, together with the fact that non-trivial type implies finite cotype, states
in other words that if $E$ has only infinite cotype, then both
$E$ and $E^*$ have only trivial type.

Evidently, any infinite dimensional Hilbert space contains $\ell_N^2$'s $1$-uniformly.
The following theorem is a variant of
Dvoretzky's theorem (see \cite{DJT}, Theorems 19.1 and 19.3 or the original paper by 
Dvoretzky \cite{DVORETZKY}), which says that Banach spaces satisfy almost
the same. 
The definition of K-convexity along with its fundamental properties can likewise be found in
\cite{DJT}, Chapter 13. 
For the purposes of this paper, one can think of K-convexity as a requirement for non-trivial
type.
Indeed, a Banach space is K-convex if and only if it has non-trivial type (\cite{DJT}, Theorem 13.3).
Furthermore, K-convexity is 
a self-dual property in the sense that a Banach space possesses 
it if and only if its dual does (\cite{DJT}, Corollary 13.7 and
Theorem 13.5).

\begin{theorem}
\label{dvothm}
  Every infinite dimensional Banach space contains $\ell_N^2$'s $\lambda$-uniformly for 
  any $\lambda > 1$.
  If the Banach space is also K-convex, then there exists a constant $C$ so that
  the $\lambda$-isomorphic copies of $\ell_N^2$'s can be chosen to be $C$-complemented.
\end{theorem}

We then turn to study the type of a space of operators. Suppose that $H$ and $E$ are Banach spaces. 
For $y\in E$ and $x^*\in H^*$ we write
\begin{equation*}
  (y \otimes x^*) x = \langle x,x^* \rangle y , \quad x\in H .
\end{equation*}
Clearly $y\otimes x^* \in \mathcal{L}(H,E)$ and $\| y \otimes x^* \| \leq \| y \| \| x^* \|$. We can also
embed $H^*$ and $E$ isometrically into $\mathcal{L}(H,E)$ 
by fixing respectively a unit vector $y\in E$ or a functional $x^*\in H^*$ with unit norm and writing
\begin{equation*}
  H^* \simeq y \otimes H^* := \{ y\otimes x^* : x^* \in H^* \} \subset \mathcal{L}(H,E)
\end{equation*}
and
\begin{equation*}
  E \simeq E \otimes x^* := \{ y\otimes x^* : y \in E \} \subset \mathcal{L}(H,E) .
\end{equation*}

The following result is most likely well-known but in lack of reference we give a proof:

\begin{prop}
  If $H$ and $E$ are infinite dimensional Banach spaces, then $\mathcal{L}(H,E)$ has only trivial type.
  \begin{proof}
    Suppose first that $H$ is K-convex and let $\lambda > 1$.
    By Theorem \ref{dvothm}, both $H$ and $E$ contain $\ell_N^2$'s $\lambda$-uniformly. 
    More precisely, there exist
    sequences $(H_N)_{N=1}^{\infty}$ and $(E_N)_{N=1}^{\infty}$ of subspaces of $H$ and $E$, such that
    each $H_N$ and $E_N$ is $\lambda$-isomorphic to $\ell_N^2$. Now, as $H$ is K-convex, 
    we may further assume that
    for some constant $C$, each $H_N$ is $C$-complemented in $H$ so that
    the projection $P_N$ onto $H_N$ has norm less or equal to $C$. We can then embed $\mathcal{L}(H_N,E_N)$
    in $\mathcal{L}(H,E)$ by extending an operator $T\in\mathcal{L}(H_N,E_N)$
    to $\widetilde{T} = TP_N$ so that $\| \widetilde{T} \| \leq C \| T \|$.
    Fix an $N$ and denote the isomorphisms from $H_N$ and $E_N$ to $\ell_N^2$ 
    by $\Lambda_N^H$ and $\Lambda_N^E$, respectively. Define
    \begin{equation*}
      \Lambda : \mathcal{L}(\ell_N^2,\ell_N^2) \to \mathcal{L}(H_N,E_N)
    \end{equation*}
    by $\Lambda (T) = (\Lambda_N^E)^{-1} T \Lambda_N^H$. Then
    $\Lambda^{-1} (S) = \Lambda_N^E S (\Lambda_N^H)^{-1}$ and
    \begin{equation*}
      \| \Lambda \| \| \Lambda^{-1} \| 
      \leq \| (\Lambda_N^E)^{-1} \| \| \Lambda_N^H \| \| \Lambda_N^E \| \| (\Lambda_N^H)^{-1} \|
      \leq \lambda^2 .
    \end{equation*}

    As every sequence in $\ell_N^{\infty}$ defines a (diagonal) operator in $\mathcal{L}(\ell_N^2,\ell_N^2)$ 
    with same operator norm, we have 
    $\ell_N^{\infty} \hookrightarrow \mathcal{L}(\ell_N^2,\ell_N^2)$ isometrically.
    Thus $\mathcal{L}(H,E)$ contains $\ell_N^{\infty}$'s $C\lambda^2$-uniformly and cannot then by Theorem
    \ref{thm:contunif} have finite cotype, and thus cannot have non-trivial type either.

    Suppose then, that $H$ is not K-convex. Then $H^*$ is not K-convex either,
    has only trivial type and contains $\ell_N^1$'s uniformly. But 
    $H^* \hookrightarrow \mathcal{L}(H,E)$ isometrically and so $\mathcal{L}(H,E)$ 
    has also only trivial type.
  \end{proof}
\end{prop}

In many questions of vector-valued harmonic analysis 
the uniform bound of a family of operators has to be replaced by its R-bound (originally defined
by Berkson and Gillespie in \cite{BERKSONGILLESPIE}).

\begin{definition*}
  A family $\mathcal{T}$ of operators in $\mathcal{L}(H,E)$ is said to be \emph{R-bounded} if there
  exists a constant $C$ such that for any $T_1,\ldots , T_N \in \mathcal{T}$ and any 
  $x_1,\ldots , x_N \in H$, regardless of $N$, we have
  \begin{equation*}
    \E  \Big\| \sum_{j=1}^N \varepsilon_j T_j x_j \Big\|^p 
    \leq C^p \E \Big\| \sum_{j=1}^N \varepsilon_j x_j \Big\|^p ,
  \end{equation*}
  for some $p\in [1, \infty )$. 
  The smallest such constant is denoted by $\mathcal{R}_p (\mathcal{T})$. We denote $\mathcal{R}_2$ by $\mathcal{R}$
  in short later on.
\end{definition*}

Basic properties of R-bounds can be found for instance in \cite{CLEMENT}. We wish only to remark that
by the Khintchine-Kahane inequality, the R-boundedness of a family does not depend on $p$, and the constants
$\mathcal{R}_p (\mathcal{T})$ are comparable. As a consequence of the inequality
$\mathcal{R}_p (\mathcal{T} + \mathcal{S}) \leq \mathcal{R}_p (\mathcal{T}) + \mathcal{R}_p (\mathcal{S})$
for any two families $\mathcal{T}$ and $\mathcal{S}$ of operators, 
every summable sequence of operators is also R-bounded:
\begin{equation*}
  \mathcal{R}_p \Big( \{ T_j \}_{j=1}^{\infty} \Big) \leq \sum_{j=1}^{\infty} \| T_j \| .
\end{equation*}

We will then compare R-boundedness and uniform boundedness.
Any R-bounded set is seen to be uniformly bounded: 
\begin{equation*}
  \sup_{T\in\mathcal{T}} \| T \|_{\mathcal{L}(H,E)} \leq \mathcal{R}_p (\mathcal{T})
\end{equation*}
for any $1\leq p < \infty$.

In Hilbert spaces also the converse holds. More generally, the following result is proven by
Arendt and Bu in \cite{AB} (while
the authors credit the proof to Pisier):

\begin{prop}
\label{Rboundsandtypes}
  Suppose that $H$ and $E$ are Banach spaces. The following are equivalent:
  \begin{enumerate}
    \item $H$ has cotype $2$ and $E$ has type $2$.
    \item Every uniformly bounded family of linear operators in $\mathcal{L}(H,E)$ is R-bounded.
  \end{enumerate}
\end{prop}

\begin{remark}
  It is clear from above that
  if $H$ and $E$ have cotype $2$ and type $2$, respectively, and
  if $\mathcal{X}\subset \mathcal{L}(H,E)$ is a Banach space whose norm dominates the operator norm, then
  all uniformly ($\mathcal{X}$-) bounded sets are also R-bounded.
\end{remark}

There are at least two natural ways to use R-boundedness for sets of vectors in $E$. One can fix a
functional $x^*$ with unit norm on a Banach space $H$ and use the embedding
$E \simeq E \otimes x^* \subset \mathcal{L}(H,E)$. Doing so, a set $\mathcal{S}$ of vectors in $E$ is R-bounded
if there exists a constant $C$ such that
\begin{equation*}
  \E \Big\| \sum_{j=1}^N \varepsilon_j (y_j \otimes x^*)x_j \Big\|^p 
  \leq C^p \E \Big\| \sum_{j=1}^N \varepsilon_j x_j \Big\|^p
\end{equation*}
for any choice of vectors $y_1, \ldots y_N \in \mathcal{S}$ and $x_1, \ldots , x_N \in H$.

In particular, one can choose the scalar field for $H$. As linear operators from the scalars to
$E$ are of the form $\lambda \mapsto \lambda y$ for some $y\in E$, it makes sense to call
a set $\mathcal{S}$ of vectors in $E$ R-bounded
if there exists a constant $C$ such that
\begin{equation*}
  \E \Big\| \sum_{j=1}^N \varepsilon_j \lambda_j y_j \Big\|^p
  \leq C^p \E \Big| \sum_{j=1}^N \varepsilon_j \lambda_j \Big|^p
\end{equation*}
for all vectors $y_1, \ldots , y_N$ in $\mathcal{S}$ and all scalars $\lambda_1, \ldots , \lambda_N$.
These two conditions are easily seen to be equivalent.

\section{The Rademacher maximal function}

Suppose that $H$ and $E$ are Banach spaces and that $\mathcal{X}\subset \mathcal{L}(H,E)$ is a Banach space
whose norm dominates the operator norm.
We are mostly interested in the case $\mathcal{X} \simeq E$, i.e. when $\mathcal{X} = E \otimes x^*$ for some
$x^*\in H^*$ or $H$ is the scalar field. 
Another typical choice for $\mathcal{X}$ is $\mathcal{L}(H,E)$ itself. Further, when $H$ is a Hilbert space, 
we can take the so-called $\gamma$-radonifying operators for our $\mathcal{X}$ (for the definition, see
Linde and Pietsch \cite{LINDE}, van Neerven \cite{GAMMARAD} or the book \cite{DJT} Chapter 12). 
Their natural norm is not
equivalent to the operator norm, thus giving us a non-trivial example of an interesting $\mathcal{X}$. Finally,
for Hilbert spaces $H_1$ and $H_2$ one can consider the Schatten-von Neumann classes $S_p(H_1,H_2)$ with
$1\leq p < \infty$ (see \cite{DJT} Chapter 4).

We will now set out to define the Rademacher maximal function.
Suppose that $(\Omega , \mathcal{F} , \mu )$ is a $\sigma$-finite measure space and
denote the corresponding Lebesgue-Bochner
space of $\mathcal{F}$-measurable $\mathcal{X}$-valued functions by 
$L^p(\mathcal{F} ; \mathcal{X})$ (or $L^p(\mathcal{X})$), $1\leq p \leq \infty$.
The space of 
strongly measurable functions $f$ for which
$1_A f$ is integrable for every set $A\in\mathcal{F}$ with finite measure, is denoted by 
$L_{\sigma}^1(\mathcal{F} ; \mathcal{X})$.

If $\mathcal{G}$ is a sub-$\sigma$-algebra of $\mathcal{F}$ such that $(\Omega, \mathcal{G}, \mu )$ is
$\sigma$-finite, there exists for every function $f\in L_{\sigma}^1(\mathcal{F} ; \mathcal{X})$
a \emph{conditional expectation} 
$\E (f | \mathcal{G}) \in L_{\sigma}^1 (\mathcal{G} ; \mathcal{X})$ with respect to $\mathcal{G}$
which is the (almost everywhere) unique strongly $\mathcal{G}$-measurable function satisfying
\begin{equation*}
  \int_A \E (f | \mathcal{G}) \D\mu = \int_A f \D\mu
\end{equation*}
for every $A\in\mathcal{G}$ with finite measure. The operator $\E ( \cdot | \mathcal{G} )$ is a contractive 
projection from $L^p(\mathcal{F} ; \mathcal{X})$ onto $L^p(\mathcal{G} ; \mathcal{X})$ for any $p\in [1,\infty ]$.
This follows immediately, if the vector-valued conditional expectation is constructed as the tensor extension of the
scalar-valued conditional expectation, which is a positive operator (see Stein
\cite{STEIN} for the scalar-valued case).

Conditional expectations satisfy Jensen's inequality:
If $\phi : \mathcal{X} \to \R$ is a convex function and $f\in L^1_{\sigma}(\mathcal{X})$ is such that
$\phi \circ f \in L^1_{\sigma}$, then
\begin{equation*}
  \phi \circ \E ( f | \mathcal{G} ) \leq \E ( \phi \circ f | \mathcal{G} )
\end{equation*}
for any sub-$\sigma$-algebra $\mathcal{G}$ of $\mathcal{F}$ (for which $(\Omega, \mathcal{G}, \mu )$ is
$\sigma$-finite). The proof in the case of a finite measure space can be found in \cite{HYTONENDIP}.

Suppose then that $(\mathcal{F}_j)_{j\in\Z}$ is a \emph{filtration}, that is, 
an increasing sequence of sub-$\sigma$-algebras of $\mathcal{F}$ such that each $(\Omega, \mathcal{F}_j, \mu )$ is
$\sigma$-finite. For a function $f\in L_{\sigma}^1(\mathcal{F} ; \mathcal{X})$, 
we denote the conditional expectations with respect to this filtration by
\begin{equation*}
  E_j f := \E (f | \mathcal{F}_j) , \quad j\in\Z .
\end{equation*}

The standard maximal function (with respect to $(\mathcal{F}_j)_{j\in\Z}$) is given by
\begin{equation*}
  Mf(\xi ) 
  = \sup_{j\in\Z} \| E_jf(\xi ) \| , \quad \xi\in\Omega ,
\end{equation*}
for functions $f$ in $L^1_{\sigma}(\mathcal{X})$.
The operator $f \mapsto Mf$ is known to be bounded from 
$L^p(\mathcal{X})$ to $L^p$ whenever $1<p\leq\infty$, regardless of $\mathcal{X}$.

\begin{definition*}
The Rademacher maximal function of a function $f\in L_{\sigma}^1(\mathcal{F} ; \mathcal{X})$ is defined by
\begin{equation*}
  M_Rf(\xi ) 
  = \mathcal{R} \Big( E_jf(\xi ) : j\in\Z \Big) , \quad \xi\in\Omega .
\end{equation*}
\end{definition*}

\begin{remark} Two immediate observations are listed below.
\begin{enumerate}
\item
The $\mu$-measurability of $M_Rf$ can be seen by studying it as the supremum over $N$ of the truncated versions
\begin{equation*}
  M_R^{(N)}f(\xi ) = \mathcal{R} \Big( E_jf(\xi ) : | j | \leq N \Big) , \quad \xi\in\Omega .
\end{equation*}
Indeed, every $M_R^{(N)}f$ is a composition of a strongly $\mu$-measurable function
\begin{equation*}
  \Omega \to \mathcal{X}^{2N+1} : \xi \mapsto (E_jf(\xi ))_{j=-N}^N
\end{equation*}
and a continuous function (we assumed that the norm of $\mathcal{X}$ dominates the operator norm)
\begin{equation*}
  \mathcal{X}^{2N+1} \to \R : \quad (T_j)_{j=-N}^N \mapsto \mathcal{R} \Big( T_j : |j| \leq N \Big) .
\end{equation*}

\item
By the properties of R-bounds we obtain the pointwise relation
$Mf \leq M_Rf$. If $H$ has cotype $2$ and $E$ has type $2$ 
it follows from Proposition \ref{Rboundsandtypes} (and the following remark) that $M_Rf \lesssim Mf$. 
This is the case in particular, when $H = L^q$ for $1 \leq q \leq 2$ and $E = L^p$ for $2 \leq p < \infty$
over some measure spaces. 

\end{enumerate}
\end{remark}
\begin{example}
Equip the Euclidean space $\R^n$ with the Borel $\sigma$-algebra and the Lebesgue measure.
For each integer $j$, let $\mathcal{D}_j$ denote a partition of $\R^n$ into
\emph{dyadic cubes} with edges of length $2^{-j}$. Suppose in addition, that every cube in $\mathcal{D}_j$
is a union of $2^n$ cubes in $\mathcal{D}_{j+1}$.
For instance, one can take the ``standard'' dyadic cubes $\mathcal{D}_j = \{ 2^{-j} ([0,1)^n + m) : m\in\Z^n \}$.
A filtration $(\mathcal{F}_j)_{j\in\Z}$ is then obtained by defining $\mathcal{F}_j$ as the $\sigma$-algebra
generated by $\mathcal{D}_j$.
We write $\langle f \rangle_Q$ for the average of an $\mathcal{X}$-valued function $f$
over a dyadic cube $Q$, that is
\begin{equation*}
  \langle f \rangle_Q = \frac{1}{|Q|}\int_Q f(\eta ) \D \eta .
\end{equation*}

Our maximal functions are now given by
\begin{equation*}
  Mf(\xi) = \sup_{Q\ni\xi} \| \langle f \rangle_Q \| \quad \text{and} \quad
  M_Rf(\xi) = \mathcal{R} \Big( \langle f \rangle_Q : Q\ni\xi \Big) , \quad \xi\in\R^n .
\end{equation*}
\end{example}

The Euclidean version of the 
Rademacher maximal function was originally studied by Hytönen, McIntosh and Portal 
\cite{HMP} via the identification $\mathcal{L}(\C , E) \simeq E$. 
They showed using interpolation that the $L^p$-boundedness 
of $f\mapsto M_Rf$ for one $p\in (1,\infty )$ implies boundedness for all 
$p$ in that range. They also provided an example of a space, namely $\ell^1$, for which
the Rademacher maximal operator is not bounded.


\begin{definition*}
Let $1 < p < \infty$.
A Banach space $\mathcal{X}\subset \mathcal{L}(H,E)$ is said to have
$\text{RMF}_p$ with respect to a given filtration on a given $\sigma$-finite measure space 
if the corresponding Rademacher maximal operator is bounded from $L^p(\mathcal{X})$ to $L^p$.
\end{definition*}

The smallest constant for which the boundedness holds will be called the $\text{RMF}_p$-constant for the given
filtration on the given measure space. 
When dealing with the Euclidean case, we occasionally drop the subscript $p$ and refer to the property
as $\text{RMF}$ with respect to $\R^n$.
Note that the $\text{RMF}_p$-property is inherited by 
closed subspaces. In particular, if $\mathcal{L}(H,E)$ has
$\text{RMF}_p$, then both $E$ and $H^*$ have it.



We will show that if $\mathcal{X}$ has $\text{RMF}_p$ with respect to the filtration of dyadic intervals on $[0,1)$,
then it has $\text{RMF}_p$ with respect to any filtration on any $\sigma$-finite measure space.
Supporting evidence is found in the Euclidean case:
If one restricts to the unit cube $[0,1)^n$ with the filtration of
dyadic cubes contained in $[0,1)^n$, it is not difficult to show that
$\text{RMF}_p$ with respect to this filtration on $[0,1)^n$ is equivalent to
$\text{RMF}_p$ with respect to the filtration of standard dyadic cubes on $\R^n$.

Martingales are later on used to study a weak type inequality for the maximal operator. In the Euclidean case,
a similar inequality can be proven with the aid of Calderón-Zygmund decomposition:
Suppose that $\mathcal{X} \subset \mathcal{L}(H,E)$ has $\text{RMF}_p$ with respect to the filtration of
dyadic cubes on $\R^n$
for some $p \in (1, \infty )$, i.e. that $M_R$ is bounded from $L^p(\mathcal{X})$ to $L^p$. 
Then there exists a constant $C$ such that for all $f\in L^1(\mathcal{X})$,
\begin{equation*}
 | \{ \xi\in\R^n : M_Rf(\xi) > \lambda \} | \leq \frac{C}{\lambda} \| f \|_{L^1(\mathcal{X})}
\end{equation*}
whenever $\lambda > 0$. The crucial part of the proof is to observe that $M_Ra$ vanishes outside a dyadic cube
containing the support of an atom $a$ (whose average is zero).

\section{RMF-property, type and cotype}
\label{RMFtypecotype}

We will now study what kind of restrictions the boundedness of the Rademacher maximal operator puts on the type and cotype
of the spaces involved. 

Unlike many other maximal operators, $M_R$ is not in general bounded from 
$L^{\infty}(\mathcal{L}(H,E))$ to $L^{\infty}$.
We actually have the following:
\begin{prop}
  The Rademacher maximal operator is bounded from $L^{\infty}(0,1;\mathcal{L}(H,E))$ 
  to $L^{\infty}(0,1)$ if and only if $H$ has cotype 2 and $E$ has type 2.
  \begin{proof}
    If $H$ has cotype 2 and $E$ has type 2, all the uniformly bounded sets are R-bounded and $M_R f \leq CMf$ 
    for all $f$ in $L^{\infty}(0,1 ; \mathcal{L}(H,E))$.
    Suppose on the contrary, that $H$ does not have cotype 2 or that $E$ does not have type 2 and fix a $C > 0$. 
    Now there exists a positive integer $N$
    and operators $T_1 , \ldots , T_N$ in $\mathcal{L}(H,E)$ with at most unit norm such that the R-bound of
    $\{ T_1, \ldots , T_N \}$ is greater than $C$. 
    We then construct an $L^{\infty}$-function on $[0,1)$ that obtains the
    operators $T_j$ as dyadic averages on an interval. Let us write $I_j = [0, 2^{j-N})$, $j=1,\ldots N$, so
    that $I_1 = [0,2^{1-N})$ is the smallest interval and $I_N = [0,1)$. We set $S_1 = T_1$ and
    \begin{equation*}
      S_j = 2T_j - T_{j-1} , \quad j=2,\ldots N .
    \end{equation*}
    Now $\| S_j \| \leq 3$ for all $j=1,\ldots , N$, so that if we define $f(\xi ) = S_1$ for $\xi \in I_1$ and
    $f(\xi ) = S_j$ for $\xi\in I_j\setminus I_{j-1}$, $j=2,\ldots , N$, we have
    $f\in L^{\infty}(0,1;\mathcal{L}(H,E))$.

\begin{figure}[h!]
\centering
  \setlength{\unitlength}{1bp}%
    \begin{picture}(238.11, 42.52)(0,0)
  \put(0,0){\includegraphics{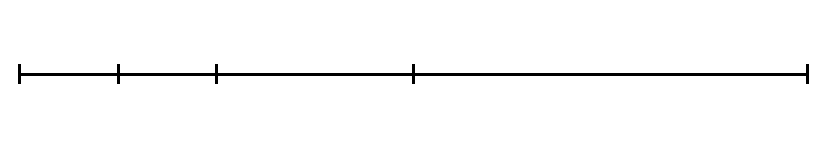}}
  \put(17.01,30.18){\fontsize{8.54}{10.24}\selectfont $S_1$}
  \put(45.35,30.18){\fontsize{8.54}{10.24}\selectfont $S_2$}
  \put(85.04,30.18){\fontsize{8.54}{10.24}\selectfont $S_3$}
  \put(164.41,30.18){\fontsize{8.54}{10.24}\selectfont $S_4$}
  \put(19.84,7.51){\fontsize{8.54}{10.24}\selectfont $I_1$}
  \put(36.85,7.51){\fontsize{8.54}{10.24}\selectfont $I_2\setminus I_1$}
  \put(79.37,7.51){\fontsize{8.54}{10.24}\selectfont $I_3\setminus I_2$}
  \put(158.74,7.51){\fontsize{8.54}{10.24}\selectfont $I_4\setminus I_3$}
  \end{picture}%
  \caption{The construction of $f$ with $N=4$}
\end{figure}
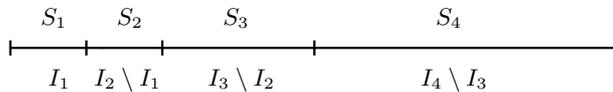

     We then look at the 
    averages of $f$ over the intervals $I_j$. 
    Obviously
    \begin{align*}
      \langle f \rangle_{I_1} &= S_1 = T_1 , \\
      \langle f \rangle_{I_2} &= \frac{S_1  + S_2}{2} = \frac{T_1 + 2T_2 - T_1}{2} = T_2  \quad \text{and} \\
      \langle f \rangle_{I_3} &= \frac{S_1 + S_2 + 2S_3}{4} = \frac{2T_2 + 4T_3 - 2T_2}{4} = T_3 .
    \end{align*}
    More generally, observing the telescopic behaviour we calculate
    \begin{equation*}
      \langle f \rangle_{I_j} = \frac{1}{2^{j-1}} \Big( S_1 + \sum_{k=1}^j 2^{k-1} S_k \Big)
      = \frac{1}{2^{j-1}} ( T_1 + 2^{j-1} T_j - T_1) = T_j ,
    \end{equation*}
    for $j=2,\ldots , N$,  
    as was desired. Thus $M_Rf > C$ on $I_1$, where $C$ was chosen arbitrarily large
    and the bound 3 for the norm of $f$ does not depend on $C$. The operator $M_R$ cannot therefore be bounded from
    $L^{\infty}(0,1;\mathcal{L}(H,E))$ to $L^{\infty}(0,1)$.

  \end{proof}
\end{prop}

Based on the counterexample from \cite{HMP} that the sequence space 
$\ell^1$ does not have RMF we prove the following statement.

\begin{prop}
\label{nontrivialtype}
  If for some $p\in (1,\infty )$, $\mathcal{L}(H,E)$ has $\text{RMF}_p$ with 
  respect to the usual dyadic filtration on $\R$, then $H$ has finite cotype and $E$ has non-trivial type.
  \begin{proof}
    Suppose on the contrary that $E$ has only trivial type. By Theorem \ref{thm:contunif} it follows that
    for some $\lambda \geq 1$ there exists a sequence $(E_N)_{N=1}^{\infty}$ of subspaces and a sequence 
    $(\Lambda_N^E)_{N=1}^{\infty}$ of isomorphisms between each $E_N$ and $\ell_N^1$ such that 
    $\| \Lambda_N^E \| \| (\Lambda_N^E)^{-1} \| \leq \lambda$.
    Let us then fix an $N$. It is shown in \cite{HMP} that there exists a function 
    $f\in L^p(0,1 ; \ell^1)$ for any
    $p \in (1,\infty )$ with the following properties:
    \begin{enumerate}
      \item $f(\xi) \in \ell_{2^N}^1$ for all $\xi\in [0,1)$,
      \item $\| f(\xi) \| = 1$ for all $\xi\in [0,1)$ so that $\| f \|_{L^p(0,1 ; \ell^1)} = 1$,
      \item $\| M_Rf \|_{L^p(0,1)} \geq C_1 \log\log N$, where the constant $C_1$ does not depend on $N$.
    \end{enumerate}
    Define then a function $g : [0,1) \to E$
    by $g(\xi) = (\Lambda_{2^N}^E)^{-1} (f(\xi ))$ and note that 
    $\| g \|_{L^p(0,1 ; E)} \leq \| (\Lambda_{2^N}^E)^{-1} \|$.
    Since $M_R$ is bounded from $L^p(0,1 ; E)$ to $L^p(0,1)$ there exists a constant $C_2$ such that
    $\| M_R g \|_{L^p(0,1)} \leq C_2 \| g \|_{L^p(0,1 ; E)}$.
    But now, since $f(\xi) = \Lambda_{2^N}^E(g(\xi))$ 
    we have
    $\| M_R f(\xi) \| \leq \| \Lambda_{2^N}^E \| \| M_R g(\xi) \|$. Thus
    \begin{equation*}
      \| M_R f \|_{L^p(0,1)} \leq \| \Lambda_{2^N}^E \| \| M_R g \|_{L^p(0,1)} 
      \leq C_2 \| \Lambda_{2^N}^E \| \| g \|_{L^p(0,1 ; E)}
      \leq C_2 \lambda
    \end{equation*}
    which gives a contradiction whenever $N$ is chosen so large that $C_1\log\log N \geq C_2\lambda$.

    The claim on finite cotype is proven similarly. Suppose on the contrary that $H$ has only infinite cotype.
    Then $H^*$ has only trivial type and one can proceed as above by defining a function
    $h : [0,1) \to H^*$ by $h(\xi ) = \Lambda_{2^N}^{H^*} (f(\xi ))$.
  \end{proof}
\end{prop}

Recall that $\mathcal{L}(H,E)$ has only trivial 
type whenever $H$ and $E$ are infinite dimensional Banach spaces.
Therefore it cannot have RMF via the identification 
$\mathcal{L}(H,E) \simeq \mathcal{L}(\C , \mathcal{L}(H,E))$.

Since $L^p$-spaces have type 2 whenever $2\leq p < \infty$, 
they have the RMF-property. We show next
that they have RMF also when $1 < p < 2$. This is implied by the following heredity property of RMF.

\begin{prop}
  Let $1 < p < \infty$.
  Suppose that $(\Sigma , \nu )$ is a $\sigma$-finite measure space and that
  $\mathcal{X} \subset \mathcal{L}(H,E)$ has $\text{RMF}_p$ with respect to 
  the usual dyadic filtration on $\R^n$. Then the space
  $L^p(\Sigma ; \mathcal{X})$ has $\text{RMF}_p$ with respect to the usual dyadic filtration on $\R^n$.
  \begin{proof}
  We use the identification
  $L^p (\R^n ; L^p(\Sigma ; \mathcal{X})) \simeq L^p (\R^n \times \Sigma ; \mathcal{X})$ and write
  \begin{equation*}
    \widetilde{M}_Rf(\xi , \eta ) = \mathcal{R} \Big( \frac{1}{|Q|} \int_Q f(\zeta , \eta ) \D\zeta : Q\ni\xi \Big) , 
    \quad (\xi , \eta ) \in \R^n \times \Sigma ,
  \end{equation*}
  for the Rademacher maximal function in the first variable. By the $\text{RMF}_p$-property of $\mathcal{X}$ we have
  for $\nu$-almost every $\eta$ that
  \begin{equation*}
  \int_{\R^n} \widetilde{M}_Rf(\xi , \eta )^p \D\xi \lesssim \int_{\R^n} \| f(\xi , \eta ) \|^p \D\xi .
  \end{equation*}
  We then calculate
  \begin{align*}
  \E \Big\| \sum_{Q\ni\xi} \varepsilon_Q \lambda_Q \langle f \rangle_Q \Big\|_{L^p(\Sigma ; \mathcal{X})}^p
  &= \int_{\Sigma} \E \Big| \sum_{Q\ni\xi} \varepsilon_Q \lambda_Q \frac{1}{|Q|}
  \int_Q f(\zeta , \eta ) \D\zeta \Big| ^p \D\nu (\eta )\\
  &\lesssim \int_{\Sigma} \widetilde{M}_Rf(\xi , \eta ) ^p \D\nu (\eta ) 
  \E \Big| \sum_{Q\ni\xi} \varepsilon_Q \lambda_Q \Big| ^p
  \end{align*}
  and so
  \begin{equation*}
  \mathcal{R} \Big( \langle f \rangle_Q : Q\ni\xi \Big) ^p 
  \lesssim \int_{\Sigma} \widetilde{M}_Rf(\xi , \eta )^p \D\nu (\eta ) .
  \end{equation*}
  Therefore,
  \begin{equation*}
  \int_{\R^n} M_Rf(\xi ) ^p \D\xi \lesssim \int_{\Sigma} \int_{\R^n}  \widetilde{M}_Rf(\xi , \eta ) ^p \D\xi \D\nu (\eta )
  \lesssim \int_{\Sigma} \int_{\R^n} \| f(\xi , \eta ) \|^p \D\xi\D\nu (\eta ) ,
\end{equation*}
so that $M_R$ is bounded from $L^p((L^p(\Sigma ; \mathcal{X}))$ to $L^p$.
\end{proof}
\end{prop}

\begin{remark}
  The previous Proposition follows also from the more general results proven in \cite{HMP}, namely that 
  both noncommutative $L^p$-spaces and all UMD function lattices have RMF.
\end{remark}


\section{Reduction to Haar filtrations}

We will show that the RMF-property is independent of the filtration and the 
underlying measure space in the following sense:

\begin{theorem}
\label{rmffilt}
  Let $1 < p < \infty$. If $\mathcal{X}$ has $\text{RMF}_p$ with respect to the filtration of dyadic intervals
  on $[0,1)$, then it has $\text{RMF}_p$ with respect to any filtration on any $\sigma$-finite
  measure space.
\end{theorem}

When this is the case, we simply say that $\mathcal{X}$ has $\text{RMF}_p$.
The proof of Theorem \ref{rmffilt} uses ideas from Maurey \cite{MAUREYUMD}, where
a similar result is proven for the UMD-property. 
We begin with the simplest possible case of filtrations of finite
algebras on finite measure spaces and proceed gradually toward more general situations. In order to do so,
we first work on measure spaces $(\Omega , \mathcal{F}, \mu )$ with $\mu (\Omega ) = 1$, that are \emph{divisible}
in the sense that any set $A\in\mathcal{F}$ with positive measure has for all 
$c\in (0,1)$ a (measurable) subset with measure $c\mu (A)$.

By a \emph{basis} of a finite subalgebra $\mathcal{G}$ of $\mathcal{F}$ 
we mean a partition of $\Omega$ into disjoint non-empty sets
$A_1, \ldots , A_m \in \mathcal{G}$ that generate the subalgebra so that each $A\in\mathcal{G}$ can be
expressed as a union of some of these $A_k$'s. Such a partition, denoted by $\text{bs}\, \mathcal{G}$, 
always exists and is unique. 
Observe that
functions measurable with respect to a finite algebra can be identified with functions defined on the basis of this
algebra (or any finer algebra). 

A filtration $(\mathcal{F}_j)_{j=1}^{\infty}$ of finite subalgebras of $\mathcal{F}$ is called a \emph{Haar filtration} if 
$\text{bs}\, \mathcal{F}_j$ consists of $j+1$ sets of positive measure. 
We also write $\mathcal{F}_0 = \{ \emptyset , \Omega \}$ so that $\text{bs}\, \mathcal{F}_0 = \{ \Omega \}$.
Furthermore,
every $\mathcal{F}_j$ is obtained from $\mathcal{F}_{j-1}$ by splitting a set $B\in\text{bs}\, \mathcal{F}_{j-1}$
into two sets $B'$ and $B''$ of positive measure. 
A Haar filtration is said to be \emph{dyadic} if
in each splitting $\mu (B')$ (and hence also $\mu (B'')$) are dyadic fractions of $\mu (B)$  and further
to be \emph{standard} if each $B$ splits into sets of equal measure.

A typical example of a filtration of finite algebras is of course the filtration of dyadic intervals on $[0,1)$.
We denote by $\mathcal{D}_j$ the finite algebra of dyadic intervals of length $2^{-j}$ on $[0,1)$ and so
\begin{equation*}
  \text{bs}\,\mathcal{D}_j 
  = \{ [ (k-1)2^{-j}, k2^{-j}) : k = 1, \ldots , 2^j \} .
\end{equation*}

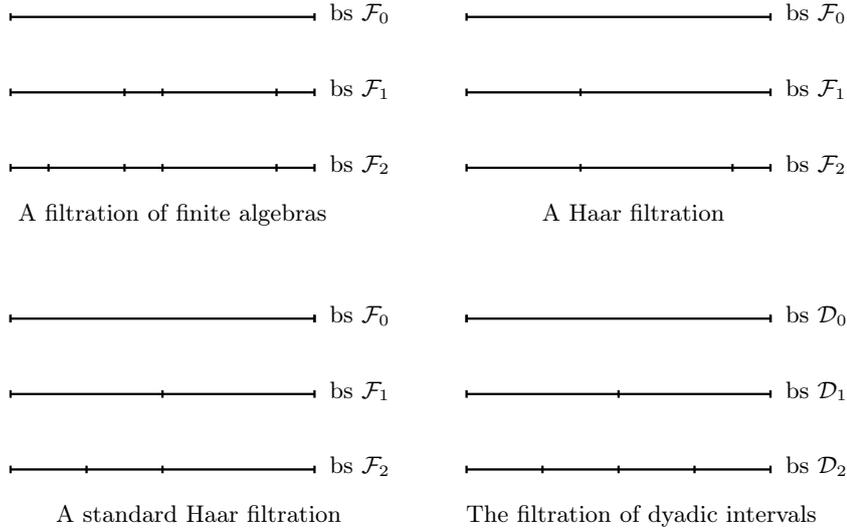
\begin{figure}[h!]
\centering
  \setlength{\unitlength}{1bp}%
  \begin{picture}(370.76, 208.35)(0,0)
  \put(0,0){\includegraphics{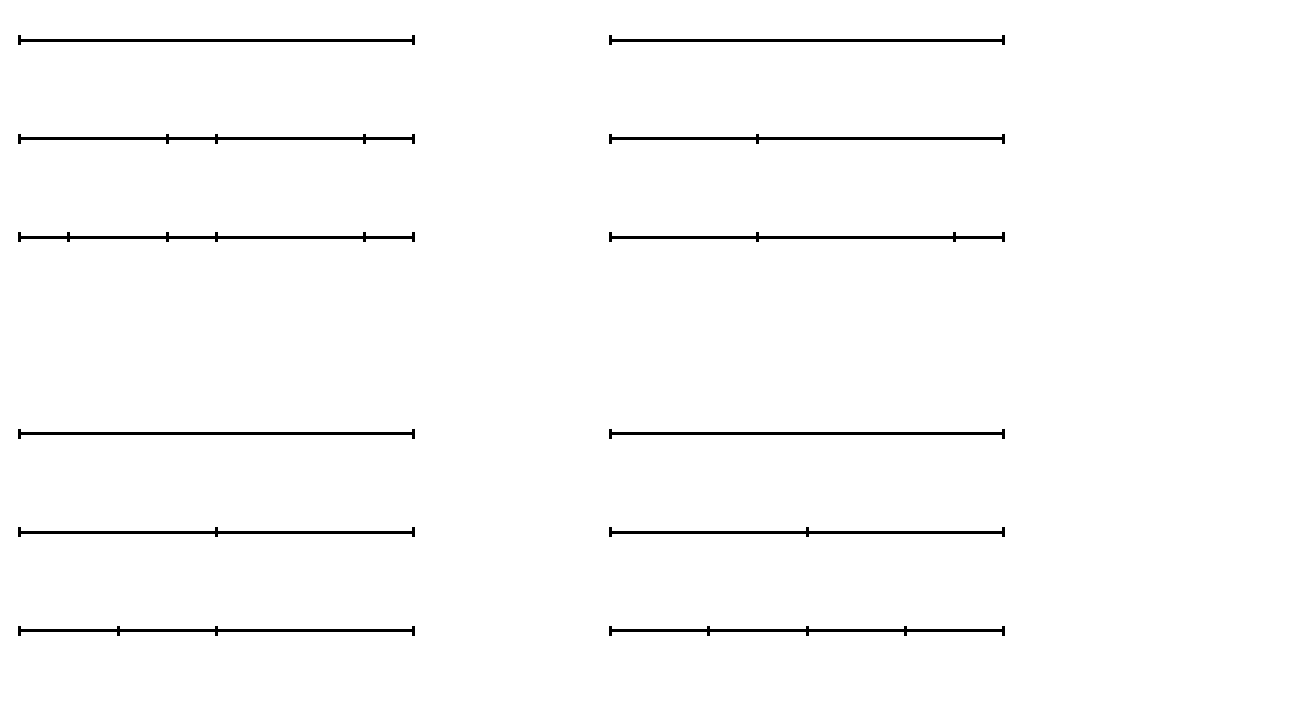}}
  \put(124.72,196.01){\fontsize{8.54}{10.24}\selectfont bs $\mathcal{F}_0$}
  \put(124.72,167.66){\fontsize{8.54}{10.24}\selectfont bs $\mathcal{F}_1$}
  \put(124.72,139.32){\fontsize{8.54}{10.24}\selectfont bs $\mathcal{F}_2$}
  \put(8.50,120.89){\fontsize{8.54}{10.24}\selectfont A filtration of finite algebras}
  \put(204.09,120.89){\fontsize{8.54}{10.24}\selectfont A Haar filtration}
  \put(294.80,196.01){\fontsize{8.54}{10.24}\selectfont bs $\mathcal{F}_0$}
  \put(294.80,167.66){\fontsize{8.54}{10.24}\selectfont bs $\mathcal{F}_1$}
  \put(294.80,139.32){\fontsize{8.54}{10.24}\selectfont bs $\mathcal{F}_2$}
  \put(124.72,82.62){\fontsize{8.54}{10.24}\selectfont bs $\mathcal{F}_0$}
  \put(124.72,54.28){\fontsize{8.54}{10.24}\selectfont bs $\mathcal{F}_1$}
  \put(124.72,25.93){\fontsize{8.54}{10.24}\selectfont bs $\mathcal{F}_2$}
  \put(294.80,82.62){\fontsize{8.54}{10.24}\selectfont bs $\mathcal{D}_0$}
  \put(294.80,54.28){\fontsize{8.54}{10.24}\selectfont bs $\mathcal{D}_1$}
  \put(294.80,25.93){\fontsize{8.54}{10.24}\selectfont bs $\mathcal{D}_2$}
  \put(22.68,7.51){\fontsize{8.54}{10.24}\selectfont A standard Haar filtration}
  \put(175.75,7.51){\fontsize{8.54}{10.24}\selectfont The filtration of dyadic intervals}
  \end{picture}%
\caption{Different filtrations of finite algebras}
\end{figure}

Suppose that $(\mathcal{F}_j)_{j=1}^N$ is a filtration of finite algebras.
By adding one set at a time (to the basis), one can construct a Haar filtration 
$(\widetilde{\mathcal{F}}_j)_{j=1}^{K_N}$ that
\begin{equation*}
      \widetilde{\mathcal{F}}_1 \subset \widetilde{\mathcal{F}}_2 \subset \cdots 
      \subset \widetilde{\mathcal{F}}_{K_1} = \mathcal{F}_1
      \subset \widetilde{\mathcal{F}}_{K_1+1} \subset \cdots \subset \widetilde{\mathcal{F}}_{K_N} = \mathcal{F}_N ,
\end{equation*}
where $K_j + 1$ is the number of sets in $\text{bs}\,\mathcal{F}_j$. Likewise, the filtration of dyadic intervals
on $[0,1)$ can be ``embedded'' in a standard Haar filtration on $[0,1)$.

Note that the $\text{RMF}_p$-constant of $\mathcal{X}$ with respect to a filtration
$(\mathcal{F}_j)_{j=1}^N$ of finite algebras
is at least the $\text{RMF}_p$-constant with respect to any ``subfiltration'' $(\mathcal{F}_{j_k})_{k=1}^M$,
where $1\leq j_{k_1} \leq \ldots \leq j_{k_M} \leq N$.
Indeed, for any $\mathcal{F}_N$-measurable $f$ we have
\begin{equation*}
  \mathcal{R} \Big( \E (f | \mathcal{F}_{j_k})(A) : 1 \leq k \leq M \Big)
  \leq \mathcal{R} \Big( \E (f | \mathcal{F}_j)(A) : 1 \leq j \leq N \Big) , \quad A\in\text{bs}\,\mathcal{F}_N ,
\end{equation*}
and the claim follows.


Two filtrations $( \mathcal{F}_j )_{j=1}^{\infty}$ and $(\widetilde{\mathcal{F}}_j )_{j=1}^{\infty}$ of finite algebras
(possibly on different measure spaces) are said to be \emph{equivalent} if there exists 
for every $j\in\Z_+$ a measure preserving bijection between 
$\text{bs}\,\mathcal{F}_j$ and $\text{bs}\,\widetilde{\mathcal{F}}_j$.
Observe that if $b$ is such a bijection from $\text{bs}\,\mathcal{F}_N$ to 
$\text{bs}\,\widetilde{\mathcal{F}}_N$, then for every $\mathcal{F}_N$-measurable $f$ we have
\begin{equation*}
  \E ( f | \mathcal{F}_j ) = \E (f\circ b^{-1} | \widetilde{\mathcal{F}}_j ) \circ b
\end{equation*}
for any $j=1,\ldots , N$.
It is a matter of calculation that the $\text{RMF}_p$-constant of $\mathcal{X}$ (if finite) is the 
same with respect to equivalent filtrations of finite algebras.

Evidently, every filtration of finite algebras on any 
measure space (of total measure one) is equivalent to a filtration
on the unit interval.
The next lemma shows that when dealing with dyadic Haar filtrations, we can choose an equivalent filtration on the
unit interval that very much resembles the filtration of dyadic intervals. The result goes back to
Maurey \cite{MAUREYUMD} and a detailed proof can be found in
Hytönen \cite{HYTONENDIP}.

\begin{lemma}
\label{boolisom}
  Every dyadic Haar filtration on any measure space 
  with total measure one is equivalent
  to a dyadic Haar filtration $(\mathcal{F}_j )_{j=1}^N$ on the unit interval such that
  $\mathcal{F}_j \subset \mathcal{D}_{K_j}$ for some integers $K_j$ and
  \begin{equation*}
    \E (f | \mathcal{F}_j ) = \E ( f | \mathcal{D}_{K_j} ), \quad 1\leq j \leq N,
  \end{equation*}
  for any $\mathcal{F}_N$-measurable $f$.
\end{lemma}

Hence, if $\mathcal{X}$ has $\text{RMF}_p$ with respect to the filtration of dyadic intervals on $[0,1)$, then it
has $\text{RMF}_p$ with respect to any dyadic Haar filtration on any measure space with total measure one and
the $\text{RMF}_p$-constant is at most the $\text{RMF}_p$-constant with respect to the filtration of dyadic
intervals.

We say that $\mathcal{X}$ has $\text{RMF}_p$ \emph{uniformly} with respect to a class of filtrations on a class
of measure spaces if the $\text{RMF}_p$-constants in question are uniformly bounded.

For the next three lemmas, fix a divisible measure space $(\Omega , \mathcal{F} , \mu )$ with $\mu (\Omega ) = 1$.
In each of the lemmas we start with a filtration $(\mathcal{F}_j)_{j=1}^{\infty}$, truncate it at a positive
integer $N$ and construct a related more ``regular'' one, whose $\sigma$-algebras we denote by
$\widetilde{\mathcal{F}}_j$. Objects corresponding to these are denoted likewise, for instance,
conditional expectations are denoted by $E_j$ and $\widetilde{E}_j$, respectively.

\begin{lemma}
  If $\mathcal{X}$ has $\text{RMF}_p$ uniformly with respect to dyadic Haar filtrations 
  on $(\Omega , \mathcal{F} , \mu )$,
  then it has $\text{RMF}_p$ uniformly with respect to all Haar filtrations on $(\Omega , \mathcal{F} , \mu )$. 
  \begin{proof}
    Suppose that $(\mathcal{F}_j)_{j=1}^{\infty}$ is a Haar filtration.
    Let $f$ be an 
    $\mathcal{F}_N$-measurable $\mathcal{X}$-valued function for a fixed
    positive integer $N$ and let $\varepsilon > 0$.
    We aim to show that 
    \begin{equation*}
      \| M_R^{(N)} f \|_{L^p} \leq r(\varepsilon , N, f) + C \| f \|_{L^p(\mathcal{X})} ,
    \end{equation*}
    where $r(\varepsilon , N, f) \to 0$ as $\varepsilon \to 0$ and $C$ depends only on 
    $\mathcal{X}$ and $p$.
    
    To construct a dyadic Haar filtration $(\widetilde{\mathcal{F}}_j)_{j=1}^N$ that approximates
    $(\mathcal{F}_j)_{j=1}^N$, we proceed inductively.
    Assume that we have constructed $\widetilde{\mathcal{F}}_{j-1}$ in
    our desired dyadic Haar filtration so that $\mu (B\Delta \widetilde{B}) < 2^{j-1-N}\varepsilon$
    whenever $\widetilde{B}$ in
    $\text{bs}\,\widetilde{\mathcal{F}}_{j-1}$ corresponds to a $B$ in $\text{bs}\,\mathcal{F}_{j-1}$.
    If $B$ in $\text{bs}\,\mathcal{F}_{j-1}$ splits into $B'$ and $B''$ in $\text{bs}\,\mathcal{F}_j$,
    then using divisibility we choose $\widetilde{B}'\subset \widetilde{B}$ whose measure
    is a dyadic fraction of $\mu (\widetilde{B})$ and which contains $\widetilde{B} \cap B'$
    while $\mu (\widetilde{B}' \setminus B') < 2^{j-1-N}\varepsilon$. Now, as
    $B'\setminus \widetilde{B}' = B'\setminus \widetilde{B} \subset B\setminus \widetilde{B}$,
    we see that
    \begin{equation*}
      \mu (B'\Delta \widetilde{B}') \leq \mu (\widetilde{B}' \setminus B')
      + \mu (B\setminus \widetilde{B} ) < \frac{\varepsilon}{2^{N-j}} .
    \end{equation*}
    Also for $\widetilde{B}'' = \widetilde{B} \setminus \widetilde{B}'$ it holds that
    \begin{equation*}
      \mu (B''\Delta \widetilde{B}'' ) 
      < \frac{\varepsilon}{2^{N-j}} .
    \end{equation*}
    We have now constructed a dyadic Haar filtration $(\widetilde{\mathcal{F}}_j)_{j=1}^N$
    for which $\mu (B\Delta \widetilde{B}) < \varepsilon$ whenever $\widetilde{B}$ corresponds
    to a $B$ in some $\mathcal{F}_j$.
    Now
    \begin{align*}
      \| M_R^{(N)} f \|_{L^p}
      &= \left( \int_{\Omega} 
      \mathcal{R} \Big( E_jf(\xi ) : 1\leq j \leq N \Big) ^p \D\mu (\xi ) \right) ^{1/p} \\
      &\leq \left( \int_{\Omega} \mathcal{R} \Big( E_jf(\xi )
      - \widetilde{E}_jf (\xi ) : 1 \leq j\leq N \Big) ^p \D\mu (\xi ) \right) ^{1/p}
      + \| \widetilde{M}_R f \|_{L^p} ,
    \end{align*}
    where the maximal operator $\widetilde{M}_R$ satisfies by assumption
    $\| \widetilde{M}_R f \|_{L^p} \leq C \| f \|_{L^p(\mathcal{X})}$ for a constant $C$ independent
    of the filtration $(\widetilde{\mathcal{F}}_j)_{j=1}^N$. 
    Estimating the R-bound in the first term by summing the norms we get
    \begin{align*}
      \left( \int_{\Omega} \mathcal{R} \Big( E_jf(\xi )
      - \widetilde{E}_jf(\xi ) : 1\leq j\leq N \Big) ^p \D\mu (\xi ) \right) ^{1/p}
      &\leq \left( \int_{\Omega} \Big( \sum_{j=1}^N \| E_jf(\xi )
      - \widetilde{E}_jf(\xi ) \| \Big) ^p \D\mu (\xi ) \right)^{1/p} \\
      &\leq \sum_{j=1}^N \| E_jf - \widetilde{E}_jf \|_{L^p(\mathcal{X})} .
    \end{align*}
    
    To estimate $\| E_jf - \widetilde{E}_jf \|_{L^p(\mathcal{X})}$ for a
    fixed $j$, recall that $\mu (B \Delta \widetilde{B}) < \varepsilon$ when
    $\widetilde{B}$ in $\text{bs}\,\widetilde{\mathcal{F}}_j$ corresponds to a
    $B$ in $\text{bs}\,\mathcal{F}_j$.
    Decomposing $\Omega$ as
    \begin{equation*}
      \Omega = \bigcup_{B\in\text{bs}\,\mathcal{F}_j} B
      = \bigcup_{B\in\text{bs}\,\mathcal{F}_j} 
      \Big( ( B \cap \widetilde{B} ) \cup ( B \setminus \widetilde{B} ) \Big) 
    \end{equation*}
    gives us
    \begin{equation*}
      \| E_jf - \widetilde{E}_jf \|_{L^p(\mathcal{X})}
      \leq \sum_{B\in\text{bs}\,\mathcal{F}_j} \Big( 
         \| 1_{B \cap \widetilde{B}} ( E_jf - \widetilde{E}_jf ) \|_{L^p(\mathcal{X})}
         + \| 1_{B \setminus \widetilde{B}} ( E_jf - \widetilde{E}_jf ) \|_{L^p(\mathcal{X})} \Big) .
    \end{equation*}
    For $\xi \in B \cap \widetilde{B}$ we have
    \begin{equation*}
      E_jf(\xi ) = \frac{1}{\mu (B )} \int_{B} f \D\mu
      \quad \text{and} \quad
      \widetilde{E}_jf(\xi ) = \frac{1}{\mu (\widetilde{B} )} \int_{\widetilde{B}} f \D\mu 
    \end{equation*}
    and thus
      using Hölder's inequality we see that
    \begin{align*}
      \| 1_{B \cap \widetilde{B}} ( E_jf - \widetilde{E}_jf ) \|_{L^p(\mathcal{X})}
      &= \mu ( B \cap \widetilde{B} ) ^{1/p}
      \Big\| \frac{1}{\mu (B )} \int_{B} f \D\mu
        - \frac{1}{\mu (\widetilde{B} )} \int_{\widetilde{B}} f \D\mu \Big\| \\
        &\leq \frac{1}{\mu (B )} \Big\| \int_{B} f \D\mu - \int_{\widetilde{B}} f \D\mu \Big\|
        + \frac{|\mu (B ) - \mu (\widetilde{B})|}{\mu (B )\mu (\widetilde{B} )}
        \Big\| \int_{\widetilde{B}} f \D\mu \Big\| \\
        &\leq \frac{1}{\mu (B )} \int_{B \Delta \widetilde{B}} \| f (\xi ) \| \D\mu (\xi )
        + \frac{|\mu (B ) - \mu (\widetilde{B})|}{\mu (B )\mu (\widetilde{B} )}
        \int_{\Omega} \| f(\xi ) \| \D\mu (\xi )\\
        &\leq \frac{\mu (B \Delta \widetilde{B})^{1-1/p}}{\mu (B)} \| f \|_{L^p(\mathcal{X})}
        + \frac{\mu (B \Delta \widetilde{B})}{\mu (B )\mu (\widetilde{B} )} 
        \| f \|_{L^p(\mathcal{X})} \\
        &\leq \Big( \varepsilon ^{1-1/p} + \frac{\varepsilon}{\mu (\widetilde{B})} \Big)
        \frac{\| f \|_{L^p(\mathcal{X})}}{\mu (B)} \\
        &\leq \Big( \varepsilon ^{1-1/p} + \frac{\varepsilon}{M - \varepsilon } \Big)
        \frac{\| f \|_{L^p(\mathcal{X})}}{M} 
      \end{align*}
     whenever $\varepsilon < M$, where 
     $M = \min \{ \mu (B) : B\in\text{bs}\,\mathcal{F}_N \}$ so that
     $\mu (B) \geq M$ and 
     $\mu (\widetilde{B}) \geq M - \varepsilon$ for each $B$ in $\text{bs}\,\mathcal{F}_j$.
    On the other hand,
    \begin{equation*}
      \| 1_{B \setminus \widetilde{B}} 
         ( E_jf - \widetilde{E}_jf ) \|_{L^p(\mathcal{X})} 
       \leq \mu (B \Delta \widetilde{B} ) ^{1/p} 
       \| E_jf - \widetilde{E}_jf \|_{L^{\infty}(\mathcal{X})}
      \leq \varepsilon ^{1/p} 2 \| f \|_{L^{\infty}(\mathcal{X})}
    \end{equation*}
    for each $B$ in $\text{bs}\,\mathcal{F}_j$.

    
    All in all, as every $\text{bs}\,\mathcal{F}_j$ contains at most
    $N+1$ sets, we have established that
    \begin{align*}
      \sum_{j=1}^N \| E_jf - \widetilde{E}_jf \|_{L^p(\mathcal{X})}
      &\leq N(N+1) \Big( \varepsilon ^{1-1/p} + \frac{\varepsilon}{M - \varepsilon } \Big) 
        \frac{\| f \|_{L^p(\mathcal{X})}}{M} 
        + N(N+1)\varepsilon^{1/p} 2 \| f \|_{L^{\infty}(\mathcal{X})} \\
      &= r(\varepsilon , N, f)
    \end{align*}
    and clearly $r(\varepsilon , N, f) \to 0$ as $\varepsilon \to 0$.

  \end{proof}
\end{lemma}

\begin{lemma}
  If $\mathcal{X}$ has $\text{RMF}_p$ uniformly with respect to Haar filtrations on $(\Omega , \mathcal{F} , \mu )$, 
  then it has $\text{RMF}_p$ uniformly with respect to filtrations of finite algebras on $(\Omega , \mathcal{F} , \mu )$.
  \begin{proof}
  This follows immediately from our earlier observations: 
  Given a filtration $(\mathcal{F}_j)_{j=1}^{\infty}$ of finite algebras and a positive
  integer $N$, we can construct a Haar filtration
  $(\widetilde{\mathcal{F}}_j)_{j=1}^{K_N}$ so that
  \begin{equation*}
    \widetilde{\mathcal{F}}_1 \subset \widetilde{\mathcal{F}}_2 \subset \cdots 
    \subset \widetilde{\mathcal{F}}_{K_1} = \mathcal{F}_1
    \subset \widetilde{\mathcal{F}}_{K_1+1} \subset \cdots \subset \widetilde{\mathcal{F}}_{K_N} = \mathcal{F}_N .
  \end{equation*}
    For any $\mathcal{F}_N$-measurable $f$ we have
    \begin{equation*}
      \mathcal{R} \Big( E_jf(A) : 1\leq j\leq N \Big)
      \leq \mathcal{R} \Big( \widetilde{E}_jf(A) : 1\leq j\leq K_N \Big) , \quad A\in\text{bs}\,\mathcal{F}_N ,
    \end{equation*}
    and the claim follows.
  \end{proof}
\end{lemma}

\begin{lemma}
\label{finitealgs}
  If $\mathcal{X}$ has $\text{RMF}_p$ uniformly with respect to filtrations of finite algebras
  on $(\Omega , \mathcal{F} , \mu )$, then it has 
  $\text{RMF}_p$ uniformly with respect to all filtrations on $(\Omega , \mathcal{F} , \mu )$.
\begin{proof}
    Suppose that $(\mathcal{F}_j )_{j=1}^{\infty}$ is a filtration, $N$ a positive integer, 
    $f$ a function in $L^p(\mathcal{F}_N ; \mathcal{X})$ and
    that $\varepsilon > 0$. We begin by choosing simple functions 
    $s_j \in L^p(\mathcal{F}_j ; \mathcal{X})$, $j=1,\ldots ,N$,
    so that
    \begin{equation*}
      \| E_jf - s_j \|_{L^p(\mathcal{X})} < \frac{\varepsilon}{2^{j+1}} .
    \end{equation*}
    For $j=1,\ldots ,N$, let $\widetilde{\mathcal{F}}_j$ be the finite algebra generated by 
    $s_1,\ldots , s_j$ and observe that
    $\widetilde{\mathcal{F}}_j \subset \widetilde{\mathcal{F}}_{j+1}$, i.e. that 
    $(\widetilde{\mathcal{F}}_j)_{j=1}^N$ is a filtration.
    Now
    \begin{align*}
      \| M_R^{(N)} f \|_{L^p}
      &= \left( \int_{\Omega} 
      \mathcal{R} \Big( E_jf(\xi ) : 1\leq j\leq N \Big) ^p \D\mu (\xi ) \right) ^{1/p} \\
      &\leq \left( \int_{\Omega} \mathcal{R} \Big( E_jf(\xi )
      - \widetilde{E}_jf(\xi ) : 1\leq j\leq N \Big) ^p \D\mu (\xi ) \right) ^{1/p}
      + \| \widetilde{M}_R f \|_{L^p} ,
    \end{align*}
    where the maximal operator $\widetilde{M}_R$ 
    satisfies $\| \widetilde{M}_R f \|_{L^p} \leq C \| f \|_{L^p(\mathcal{X})}$ for a constant $C$ independent
    of the filtration $(\widetilde{\mathcal{F}}_j)_{j=1}^N$. This independence is crucial, as 
    $\widetilde{\mathcal{F}}_j$'s arose from $f$.

    We then estimate
    \begin{align*}
      \left( \int_{\Omega} \mathcal{R} \Big( E_jf(\xi )
      - \widetilde{E}_jf(\xi ) : 1\leq j\leq N \Big) ^p \D\mu (\xi ) \right) ^{1/p}
      &\leq \left( \int_{\Omega} \Big( \sum_{j=1}^N \| E_jf(\xi )
      - \widetilde{E}_jf(\xi ) \| \Big) ^p \D\mu (\xi ) \right)^{1/p} \\
      &\leq \sum_{j=1}^N \| E_jf - \widetilde{E}_jf \|_{L^p(\mathcal{X})} \\
      &\leq \sum_{j=1}^N \Big( \| E_jf - s_j \|_{L^p(\mathcal{X})} 
      + \| \widetilde{E}_jf - s_j \|_{L^p(\mathcal{X})} \Big) .     
    \end{align*}
    Furthermore, since 
    \begin{equation*}
      \| \widetilde{E}_jf - s_j \|_{L^p(\mathcal{X})} = \| \widetilde{E}_jf - \widetilde{E}_js_j \|_{L^p(\mathcal{X})}
      = \| \widetilde{E}_j (E_jf 
      - s_j) \|_{L^p(\mathcal{X})} 
      \leq \| E_jf - s_j \|_{L^p(\mathcal{X})} 
    \end{equation*}
    we get
    \begin{equation*}
      \left( \int_{\Omega} \mathcal{R} \Big( E_jf(\xi )
      - \widetilde{E}_jf(\xi ) : 1\leq j\leq N \Big) ^p \D\mu (\xi ) \right) ^{1/p}
      \leq 2 \sum_{j=1}^N \| E_jf - s_j \|_{L^p(\mathcal{X})} 
      < \sum_{j=1}^N \frac{\varepsilon}{2^j} < \varepsilon .
    \end{equation*}
  \end{proof}
\end{lemma}

We then show that the assumption on divisibility can be dropped.

\begin{lemma}
  If $\mathcal{X}$ has $\text{RMF}_p$ with respect to any filtration on any divisible 
  measure space with total measure one, 
  then it has $\text{RMF}_p$ with respect to any filtration on any measure space with total measure one.

  \begin{proof}
    Suppose that $(\mathcal{F}_j)_{j=1}^{\infty}$ is a filtration on a not necessarily divisible measure space 
    $(\Omega , \mathcal{F}, \mu )$ with $\mu (\Omega ) = 1$.
    Now the $\sigma$-algebras $\widetilde{\mathcal{F}}_j = \{ F\times [0,1] : F\in\mathcal{F}_j \}$ form a filtration on 
    the product of $(\Omega , \mathcal{F}, \mu )$ and the unit interval with Lebesgue measure, which obviously
    constitutes a divisible measure space.
    For a function
    $f\in L^p(\Omega ; \mathcal{X})$ we put $\widetilde{f}(\xi , t) = f(\xi ), (\xi , t)\in \Omega\times [0,1]$,
    and observe that $\| \widetilde{f} \|_{L^p(\mathcal{X})} = \| f \|_{L^p(\mathcal{X})}$.
    Also $\widetilde{E}_j \widetilde{f} (\xi , t) = E_j f (\xi )$ 
    for all $(\xi , t)\in\Omega \times [0,1]$, and so 
    $\| \widetilde{M}_R \widetilde{f} \|_{L^p(\mathcal{X})} = \| M_R f \|_{L^p(\mathcal{X})}$.
  \end{proof}
\end{lemma}

The results follow immediately for finite measure spaces. Suppose that $(\Omega , \mathcal{F} , \mu )$ is such. Then
the above argument applies to the measure
$\mu (\Omega )^{-1} \mu$ on $(\Omega , \mathcal{F} )$ and evidently the conditional 
expectations are the same in these two measure spaces. 
Thus the Rademacher maximal operator remains unaltered and the inequality stating the
boundedness is only a matter of scaling by $\mu (\Omega )^{-1}$.

Suppose then that $\mathcal{X}$ has $\text{RMF}_p$ uniformly with respect to any filtration on any 
finite measure space and let
$(\Omega , \mathcal{F} , \mu )$ be a $\sigma$-finite measure space with a filtration $(\mathcal{F}_j ) _{j=1}^{\infty}$.
Since $\mathcal{F}_1$ is $\sigma$-finite (by definition), 
we can write $\Omega$ as a union of disjoint sets $A_k\in\mathcal{F}_1$, 
$k\in\Z_+$, each with finite measure. 
Let us define for positive integers $k$ the finite measures $\mu_k (A) = \mu (A \cap A_k)$ on $\mathcal{F}$.
The conditional expectation of a function $f\in L^p(\Omega ; \mathcal{X})$ with
respect to $\mathcal{F}_j$ and $\mu_k$ is simply the conditional expectation of $1_{A_k}f$ with respect to
$\mathcal{F}_j$ which further equals $1_{A_k} E_j f$. In symbols
\begin{equation*}
  E_j^{(k)} f = 1_{A_k} E_jf ,
\end{equation*}
where $E_j^{(k)} f$ denotes the conditional expectation of $f$ with respect to $\mathcal{F}_j$ and $\mu_k$.
Thus
\begin{align*}
  \| M_Rf \|_{L^p}^p 
  &= \sum_{k=1}^{\infty} \int_{A_k} 
  \mathcal{R}\Big( E_j f(\xi ) : j\in\Z_+ \Big) ^p \D \mu(\xi ) \\
  &= \sum_{k=1}^{\infty} \int_{A_k}  
  \mathcal{R}\Big( E_j^{(k)} f(\xi ) : j\in\Z_+ \Big) ^p \D \mu_k (\xi ) \\
  &\leq \sum_{k=1}^{\infty} C^p \int_{A_k} \| f(\xi ) \|^p \D\mu_k (\xi ) \\
  &= C^p \| f \|_{L^p(\mathcal{X} )}^p .
\end{align*}

So far we have only considered filtrations indexed by positive integers. Suppose that $\mathcal{X}$ has $\text{RMF}_p$ with
respect to any filtration indexed by $\Z_+$ on any $\sigma$-finite measure space and let
$(\mathcal{F}_j)_{j\in\Z}$ be a filtration on $(\Omega , \mathcal{F}, \mu )$. Then $\mathcal{X}$ has 
$\text{RMF}_p$ with respect to $(\mathcal{F}_j)_{j=-N}^{\infty}$ with a constant
independent of $N$ and thus by monotone convergence theorem with respect to $(\mathcal{F}_j)_{j\in\Z}$.

This concludes the proof of Theorem \ref{rmffilt}.

\section{The weak RMF-property}

We start by recalling some terminology.
A \emph{stochastic process} (a sequence of random variables on some probability space)
$X = (X_j)_{j=1}^{\infty}$ is always \emph{adapted} to the
filtration $(\mathcal{F}_j)_{j=1}^{\infty}$, where $\mathcal{F}_j$ is the $\sigma$-algebra 
$\sigma (X_1, \ldots , X_j)$ generated by $X_1, \ldots , X_j$, in the sense
that each $X_j$ is $\mathcal{F}_j$-measurable. We call a sequence of $L^1$-random variables
a \emph{martingale} if $\E (X_k | \mathcal{F}_j) = X_j$ whenever $j\leq k$.

Note that for any martingale $X = (X_j)_{j=1}^{\infty}$ we have $\E X_j = \E X_k$ for all $j,k\in\Z_+$.
It is customary to write $\mathcal{F}_0$ for the trivial $\sigma$-algebra and $X_0$ for the common expectation of $X_j$'s.
By defining $Y_j = X_j - X_0$ one can restrict to martingales $Y = (Y_j)_{j=1}^{\infty}$ for which
$Y_0 = \E Y_j = 0$.


We say that a stochastic process $X = (X_j)_{j=1}^{\infty}$ is \emph{$L^p$-bounded} for $p\in [1,\infty )$ if
$\| X \|_p^p := \sup_{j\in\Z_+} \E \| X_j \|^p < \infty$ and for $p=\infty$ if the infimum
$\| X \|_{\infty}$ of all $C$ for which every $\| X_j \| \leq C$ almost surely, is finite.
A stochastic process $X = (X_j)_{j=1}^{\infty}$ is said to be \emph{simple} if the algebras
$\mathcal{F}_j$ are finite (i.e. if the random variables $X_j$ are simple). A simple martingale is called
a \emph{(dyadic/standard) Haar martingale} if the algebras $\mathcal{F}_j$ form a 
(dyadic/standard) Haar filtration.


Given a martingale $(X_j)_{j=1}^N$ we define its \emph{difference sequence} $(D_j)_{j=1}^N$ by
$D_j = X_j - X_{j-1}$ for $j\geq 1$. 
Furthermore, if
$v = (v_j)_{j=1}^{\infty}$ is a real $L^{\infty}$-bounded stochastic process (on the same probability space),
we define 
\begin{equation*}
  (v \star X)_j = \sum_{k=1}^j v_k D_k, \quad j \in\Z_+ .
\end{equation*}
If $v$ is \emph{predictable} with respect to $X$ in the
sense that each $v_j$ is $\mathcal{F}_{j-1}$-measurable 
(and $v_1$ is constant almost surely),
then the \emph{martingale transform}
$v\star X = ((v\star X)_j)_{j=1}^{\infty}$ is itself a martingale.

\begin{definition*}
  Let $1 < p < \infty$. A Banach space $E$ is said to have $\text{UMD}_p$ if there exists a constant $C$ such that
  for every $L^p$-martingale $X = (X_j)_{j=1}^N$ in $E$ we have
  \begin{equation*}
    \E \| (\varepsilon \star X)_N \|^p \leq C^p \E \| X_N \|^p
  \end{equation*}
  whenever $\varepsilon = (\varepsilon_j )_{j=1}^N$ is a sequence of signs $\{ -1, 1 \}$.
\end{definition*}



This property is independent of $p$ in the sense that if a Banach space has $\text{UMD}_p$ for one $p \in (1,\infty )$
then it has $\text{UMD}_p$ for all $p \in (1,\infty )$ (see Maurey \cite{MAUREYUMD}). 
Thus the parameter $p$ can be omitted from the definition. 


One can ask how the RMF-property relates to the UMD-property. First of all,
every UMD-space can be shown to be reflexive (see for instance \cite{MAUREYUMD}). 
Our typical example $\mathcal{L}(H,E)$ is usually non-reflexive, but has RMF at least when $H$ has cotype 2 and $E$
has type 2. More interestingly, James constructed in \cite{JAMES} a non-reflexive
Banach space $E$ with type 2. Thus $E \hookrightarrow \mathcal{L}(H,E)$ can have RMF without being a UMD-space.
Bourgain showed in \cite{BOURGAIN} that the Schatten-von Neumann class $S_p(H_1,H_2)$ is UMD for $1 < p < \infty$.
As $H_1$ and $H_2$ are spaces of type and cotype $2$, it follows from our earlier observations 
that $S_p(H_1,H_2)$ has RMF as a subspace of $\mathcal{L}(H_1,H_2)$. It has also been
shown in \cite{HMP} that $S_p(H_1,H_2)$ has RMF as $\mathcal{L}(\C , S_p(H_1,H_2))$.

Let $\mathcal{X} \subset \mathcal{L}(H,E)$ be a Banach space whose norm dominates the operator norm.
For a stochastic process $X = (X_j)_{j=1}^{\infty}$ in $\mathcal{X}$ we define the 
Doob and Rademacher maximal functions by
\begin{equation*}
  X^* = \sup_{j\in\Z_+} \| X_j \|  \quad \quad \text{and} \quad \quad
  X_R^* = \mathcal{R} \Big( X_j : j\in\Z_+ \Big) ,
\end{equation*}
respectively.

The boundedness properties of Doob's maximal function are well-known:
Every $L^p$-bounded martingale $X$ satisfies
\begin{equation*}
  \E | X^* |^p \leq (p')^p \| X \|_p^p ,
\end{equation*}
where $p'$ is the Hölder conjugate of $p$ and $1 < p < \infty$. 
Furthermore, for every $L^1$-bounded martingale $X$ we have
\begin{equation*}
  \prob (X^* > \lambda ) \leq \frac{1}{\lambda} \| X \|_1
\end{equation*}
whenever $\lambda > 0$.

Recall that the $\text{RMF}_p$-property is independent of the filtration and of the underlying measure space
in the sense of the previous section (Theorem \ref{rmffilt}). 
Regarding the unit interval as a probability space on which the conditional
expectations with respect to dyadic intervals define martingales, we see that
$\mathcal{X}$ has $\text{RMF}_p$ if and only if there exists a constant $C$ such that
\begin{equation*}
    \E | X_R^* |^p \leq C^p \| X \|_p^p
  \end{equation*}
for any $L^p$-bounded martingale $X$ in $\mathcal{X}$.

Applying ideas from Burkholder \cite{BURKHOLDERUMD} we will show that
$\mathcal{X}$ has $\text{RMF}_p$ for some $p\in (1,\infty )$ if and only if it has \emph{weak RMF} i.e. if
there exists a constant $C$ such that all $L^1$-bounded martingales $X$ in $\mathcal{X}$ satisfy
\begin{equation}
\label{weaktype}
  \prob (X_R^* > \lambda ) \leq \frac{C}{\lambda} \| X \|_1
\end{equation}
whenever $\lambda > 0$.

To show the necessity of the weak type inequality (\ref{weaktype}) we invoke the Gundy decomposition 
(see Gundy \cite{GUNDY} for the original proof). 

\begin{theorem} (Gundy decomposition)
  Suppose that $X$ is an $L^1$-bounded martingale in $\mathcal{X}$ and that $\lambda > 0$.
  There exists a decomposition $X=G+H+B$ of $X$ into martingales $G$, $H$ and $B$ which satisfy
  \begin{enumerate}
    \item $\| G \|_1 \leq 4 \| X \|_1$ \quad and \quad $\| G \|_{\infty} \leq 2\lambda$,
    \item $\E \| H_1 \| + \sum_{j=2}^{\infty} \E \| H_j - H_{j-1} \| \leq 4 \| X \|_1$, \quad 
          $(H = (H_j)_{j=1}^{\infty})$,
    \item $\prob ( B^* > 0 ) \leq \frac{3}{\lambda} \| X \|_1$ .
  \end{enumerate}
\end{theorem}


\begin{prop}
\label{gundyweak}
   If $\mathcal{X}$ has $\text{RMF}_p$ for some $p \in (1,\infty )$, then it has weak RMF.
  \begin{proof}
    Taking the Gundy decomposition of $X$ at height $\lambda$ we may write
    \begin{equation*}
      \prob ( X_R^* > \lambda ) \leq \prob ( B_R^* > \lambda / 3 ) + \prob ( H_R^* > \lambda / 3 ) 
      + \prob ( G_R^* > \lambda / 3 ) ,
    \end{equation*}
    and estimate each term separately. Firstly $\prob ( B_R^* > 0 ) = \prob ( B^* > 0 )$, since
    $B_R^* = 0$ if and only if $B^* = 0$. Thus
    \begin{equation*}
      \prob ( B_R^* > \lambda / 3 ) \leq \prob ( B_R^* > 0 )  = \prob ( B^* > 0 ) \leq \frac{3}{\lambda} \| X \|_1 .
    \end{equation*}
    Secondly, 
    \begin{equation*}
      H_R^* = \mathcal{R} \Big(  H_j : j\in\Z_+ \Big) 
      = \mathcal{R}\Big( \sum_{k=1}^j (H_k - H_{k-1}) : j\in\Z_+ \Big) 
      \leq \sum_{j=1}^{\infty} \| H_j - H_{j-1} \| ,
    \end{equation*}
    where the last inequality follows from a simple rearrangement of sums.
    Hence
    \begin{align*}
      \prob ( H_R^* > \lambda / 3 ) &\leq \prob 
      \Big( \sum_{j=1}^{\infty} \| H_j - H_{j-1} \| > \frac{\lambda}{3} \Big)\\
      &\leq \frac{3}{\lambda} 
      \E \sum_{j=1}^{\infty} \| H_j - H_{j-1} \| \\
      &= \frac{3}{\lambda} \sum_{j=1}^{\infty} 
      \E \| H_j - H_{j-1} \| 
      \leq \frac{12}{\lambda} \| X \|_1 .
    \end{align*}
    Thirdly,
    \begin{equation*}
      \prob ( G_R^* > \lambda / 3 ) \leq \left( \frac{3}{\lambda} \right)^p \E | G_R^* |^p 
      \leq C \left( \frac{3}{\lambda} \right)^p \| G \|_p^p 
      \leq C \frac{3^p 2^{p-1}}{\lambda} \| G \|_1 
      \leq C \frac{3^p 2^{p+1}}{\lambda} \| X \|_1 ,
    \end{equation*}
    where the property $\| G \|_{\infty} \leq 2\lambda$ was used to deduce that
    \begin{equation*}
      \| G \|_p^p 
      = \sup_{j\in\Z_+} \E \| G_j \|^p \\
      \leq \| G \|_{\infty}^{p-1} \sup_{j\in\Z_+} \E \| G_j \| \\
      \leq (2\lambda )^{p-1} \| G \|_1 .
    \end{equation*}
  \end{proof}
\end{prop}

We then turn to the converse. We obtain the desired results for standard Haar martingales, but recalling
the earlier reduction, this will not be a restriction.
The argument is based on a ``good-$\lambda$ inequality'' (Lemma \ref{goodlambda}) which
says roughly that the chance of $X_R^*$ being large while $X^*$ diminishes is vanishingly small.

\begin{lemma}
\label{predictable}
  If $X = (X_j)_{j=1}^{\infty}$ is a standard Haar martingale, then
  $( \| D_j \|)_{j=1}^{\infty}$ is predictable with respect to $X$.
  \begin{proof}
    For every $j \geq 1$ there is exactly one event $B\in \text{bs}\mathcal{F}_{j-1}$ on which
    $X_j - X_{j-1}$ is non-zero. As $B = B_1 \cup B_2$ for some $B_1, B_2 \in \text{bs}\mathcal{F}_j$
    with $\prob (B_1) = \prob (B_2)$ and $\E(X_j - X_{j-1} | \mathcal{F}_{j-1}) = 0$, 
    there exists a $T\in\mathcal{X}$ such that
    $X_j - X_{j-1} = 1_{B_1}T - 1_{B_2}T$. Consequently,
    \begin{equation*}
     \| D_j \| = \| X_j - X_{j-1} \| = 1_{B_1} \| T \| + 1_{B_2} \| T \| = 1_B \| T \|
    \end{equation*}
    and so $\| D_j \|$ is $\mathcal{F}_{j-1}$-measurable.
  \end{proof}
\end{lemma}

We will need the concept of a stopping time:
We say that a random variable $\tau$ in $\Z_+ \cup \{ \infty \}$ is a \emph{stopping time} with respect to a stochastic
process $X$ if $\{ \tau = j \}$ is in $\mathcal{F}_j$ for every positive integer $j$. In this case we define
\begin{equation*}
  X_{\tau} = \sum_{j=1}^{\infty} 1_{\{ \tau = j \} } X_j .
\end{equation*}
Observe that $X_{\tau} = 0$ when $\tau = \infty$. An easy calculation shows that if $\tau$ is a stopping time
with respecto to an $L^1$-bounded martingale $X$, then $\E \| X_{\tau} \| \leq \| X \|_1$.

\begin{lemma}
\label{goodlambda}
  Suppose that $\mathcal{X}$ has weak RMF.
  Then for all $\delta \in (0,1)$ and $\beta > 2\delta +1$ 
  there exists an $\alpha (\delta ) > 0$ which tends to zero as $\delta \searrow 0$ and which is such that
  for all $L^p$-bounded standard Haar martingales $X$ in $\mathcal{X}$ we have
  \begin{equation*}
    \prob \Big( X_R^* > \beta\lambda , \, X^* \leq \delta \lambda \Big) \leq \alpha (\delta ) \prob ( X_R^* > \lambda ) ,
  \end{equation*}
  whenever $\lambda > 0$.
  \begin{proof}
    Let $X = (X_j)_{j=1}^{\infty}$ be an $L^p$-bounded standard Haar martingale in $\mathcal{X}$. 
    Define the stopping times
    \begin{align*}
      \tau_1 &= \min \Big\{ j\in\Z_+ : \mathcal{R}\Big( X_k : 1\leq k\leq j \Big) > \lambda \Big\} \\
      \tau_2 &= \min \Big\{ j\in\Z_+ : \mathcal{R}\Big( X_k : 1\leq k\leq j \Big) > \beta\lambda \Big\} \\
      \sigma &= \min \Big\{ j\in\Z_+ : \| X_j \| > \delta \lambda \quad \text{or} \quad
      \| D_{j+1} \| > 2\delta \lambda \Big\}
    \end{align*}
    where Lemma \ref{predictable} guarantees that $\{ \sigma = j \} \in \mathcal{F}_j$ for each
    $j\in\Z_+$.
    Define then
    \begin{equation*}
      v_j = 1_{\{ \tau_1 < j \leq \tau_2 \wedge \sigma \} } ,
    \end{equation*}
    and note that $\{ \tau_1 < j \leq \tau_2 \wedge \sigma \}$ is the intersection of
    the complements of $\{ \tau_2\wedge\sigma > j \}$ and $\{ \tau_1 > j-1 \}$, both of which lie in
    $\mathcal{F}_{j-1}$.
    Hence $v = (v_j)_{j=1}^{\infty}$ is predictable and so $v\star X$ is a martingale.
    When $\tau_1 < \tau_2 \wedge \sigma$ we calculate
    \begin{equation*}
      (v\star X)_j = \sum_{k=1}^j v_kD_k 
      = \sum_{\tau_1 < k \leq \tau_2 \wedge \sigma \wedge j} (X_k - X_{k-1})
      =
      \begin{cases}
        0, \quad 1\leq j \leq \tau_1 , \\
        X_j - X_{\tau_1}, \quad \tau_1 < j \leq \tau_2 \wedge \sigma , \\
        X_{\tau_2\wedge \sigma} - X_{\tau_1}, \quad j > \tau_2 \wedge \sigma .
      \end{cases} 
    \end{equation*}
    We first show that
    \begin{equation*}
      \{ X_R^* > \beta\lambda , \, X^* \leq \delta\lambda \}
      \subset \{ (v\star X)_R^* > (\beta - 2\delta - 1) \lambda \} .
    \end{equation*}

    
    Suppose that $X_R^* > \beta\lambda$ and $X^* \leq \delta\lambda$. Now $\tau_2 < \infty$ and as
    $\| D_{j+1} \| \leq \| X_{j+1} \| + \| X_j \| \leq 2\delta\lambda$ for all $j$, we also have
    $\sigma = \infty$.
    Since for every $j$
    \begin{equation*}
      \mathcal{R}\Big( X_k : 1\leq k\leq j \Big)
      \leq \mathcal{R}\Big( X_k : 1\leq k\leq j-1 \Big) + \| D_j \| ,
    \end{equation*}
    we have
    \begin{equation*}
      \mathcal{R}\Big( X_k : 1\leq k \leq \tau_2 - 1 \Big) 
      \geq \mathcal{R}\Big( X_k : 1\leq k \leq \tau_2 \Big) - \| D_{\tau_2} \|
      > (\beta - 2\delta ) \lambda > \lambda .
    \end{equation*}
    Thus $\tau_1 < \tau_2$ and
    \begin{equation*}
      (v\star X)_j =
      \begin{cases}
        0, \quad 1\leq j \leq \tau_1 , \\
        X_j - X_{\tau_1}, \quad \tau_1 < j \leq \tau_2 , \\
        X_{\tau_2} - X_{\tau_1}, \quad j > \tau_2 .
      \end{cases}
    \end{equation*}
    Hence
    \begin{align*}
      (v\star X)_R^* &= \mathcal{R} \Big( X_j - X_{\tau_1} : \tau_1 < j \leq \tau_2 \Big) \\
      &\geq \mathcal{R} \Big( X_j : \tau_1 < j \leq \tau_2 \Big) - \| X_{\tau_1} \| \\
      &\geq \mathcal{R} \Big( X_j : 1 \leq j \leq \tau_2 \Big) 
      - \mathcal{R} \Big( X_j : 1 \leq j \leq \tau_1 \Big) - \| X_{\tau_1} \| \\
      &\geq \mathcal{R} \Big( X_j : 1 \leq j \leq \tau_2 \Big) 
      - \mathcal{R} \Big( X_j : 1 \leq j < \tau_1 \Big) - 2\| X_{\tau_1} \| \\
      &> \beta\lambda - \lambda - 2  \delta\lambda \\
      &> (\beta - 2\delta - 1)\lambda ,
    \end{align*}
    as required.

    We then aim to find a suitable upper bound for $\| v\star X \|_1$. To do this, consider cases
    $\{ \tau_1 < \tau_2 \wedge \sigma \}$ and $\{ \tau_1 \geq \tau_2 \wedge \sigma \}$ separately. 
    In the former case, an earlier calculation gives
    \begin{equation*}
      \| (v\star X)_j \| \leq \| X_{\tau_2 \wedge \sigma \wedge j} \| + \| X_{\tau_1} \| ,
    \end{equation*}
    where $\| X_{\tau_1} \| \leq \delta\lambda$. Furthermore
    \begin{equation*}
      \| X_{\tau_2 \wedge \sigma \wedge j} \| \leq \| X_{\tau_2 \wedge \sigma \wedge j - 1} \| + 
      \| D_{\tau_2\wedge \sigma \wedge j} \| \leq \delta\lambda + 2\delta\lambda
    \end{equation*}
    and so
    $\| (v\star X)_j \| \leq 4\delta\lambda$ for all $j\in\Z_+$. 
    In the latter case each $v_j = 0$ and so $(v\star X)_j = 0$. This happens in 
    particular in the event
    $\{ \tau_1 = \infty \} = \{ X_R^* \leq \lambda \}$. Thus in conclusion
    \begin{equation*}
      (v\star X)^* \leq 4\delta\lambda 1_{\{ \tau_1 < \infty \} }
    \end{equation*}
    and so
    \begin{equation*}
      \| v\star X \|_1 \leq \E (v\star X)^* \leq 4\delta\lambda \prob (X_R^* > \lambda ) .
    \end{equation*}
    Putting all these estimates together we get
    \begin{align*}
      \prob \Big( X_R^* > \beta\lambda , \, X^* \leq \delta \lambda \Big)
      &\leq \prob \Big( (v\star X)_R^* > (\beta - 2\delta - 1) \lambda \Big) \\
      &\leq \frac{C}{(\beta - 2\delta - 1) \lambda} \| v\star X \|_1 \\
      &\leq \frac{4C\delta}{(\beta - 2\delta - 1)} \prob (X_R^* > \lambda ) .
    \end{align*}
    Fixing a $\beta > 2\delta + 1$ we may take
    \begin{equation*}
      \alpha (\delta ) = \frac{4C\delta}{(\beta - 2\delta - 1)} .
    \end{equation*}
  \end{proof}
\end{lemma}

The previous lemma allows us to deduce the strong type inequality from the weak type inequality:

\begin{prop}
\label{trick}
  Suppose that $\mathcal{X}$ has weak RMF and
  let $1 < p < \infty$. Then there exists a constant $C$ such that 
  for any $L^p$-bounded standard Haar martingale $X$ in $\mathcal{X}$ we have
  $\E |X_R^*|^p \leq C^p \| X \|_p^p$.  
  \begin{proof}
   Let $X = (X_j)_{j=1}^N$ be a standard Haar martingale in $\mathcal{X}$ (note that it suffices to prove the claim for
   finite martingales independently of $N$).
   We apply the good-$\lambda$ inequality and write
    \begin{align*}
      \E | X_R^* |^p &= \beta^p \int_0^{\infty} p \lambda ^{p-1} \prob (X_R^* > \beta\lambda ) \D\lambda \\
      &\leq \beta^p \alpha (\delta ) \int_0^{\infty} p \lambda ^{p-1} \prob ( X_R^* > \lambda ) \D \lambda
           + \beta^p \int_0^{\infty} p \lambda ^{p-1} \prob ( X^* > \delta \lambda ) \D\lambda \\
      &= \beta^p \alpha (\delta ) \E | X_R^* |^p + \frac{\beta^p}{\delta ^p} \E | X^* |^p ,
    \end{align*}
    where $\E | X^* |^p \leq C^p \| X \|_p^p$ and $\E | X_R^* |^p$ is finite. 
    Choosing $\delta$ so small that $\beta^p \alpha (\delta ) < 1$ we get
    \begin{equation*}
      \E | X_R^* |^p \leq \frac{\beta^p C^p}{(1 - \beta^p \alpha (\delta ) ) \delta ^p} \| X \|_p^p .
    \end{equation*}
  \end{proof}
\end{prop}

We collect our results as follows:

\begin{theorem}
\label{weakrmf}
The following conditions are equivalent:
\begin{enumerate}
  \item $\mathcal{X}$ has $\text{RMF}_p$ for all $p\in (1,\infty )$.
  \item $\mathcal{X}$ has $\text{RMF}_p$ for some $p\in (1,\infty )$.
  \item $\mathcal{X}$ has weak RMF.
\end{enumerate}
\begin{proof}
  Trivially the first condition implies the second. That the third follows from the second was 
  Proposition \ref{gundyweak}. In Proposition \ref{trick} we showed that the weak RMF-property implies that
  for any $p\in (1,\infty )$,
  $\E |X_R^*|^p \lesssim \| X \|_p^p$ whenever $X$ is an $L^p$-bounded standard Haar martingale in $\mathcal{X}$.
  As was noted before, the filtration of dyadic intervals on $[0,1)$ can be ``embedded'' in a standard Haar filtration.
  Thus the weak RMF-property is sufficient for the $L^p$-boundedness, $1 < p < \infty$, of the Rademacher maximal operator
  on the unit interval. By Theorem \ref{rmffilt} this implies $\text{RMF}_p$ for all $p\in (1,\infty )$.
\end{proof}
\end{theorem}

\section{RMF-property and concave functions}

The existence of a biconcave function $v : E \times E \to \R$ for which
\begin{equation*}
  v(x,y) \geq \Big\| \frac{x+y}{2} \Big\|^p - C^p \Big\| \frac{x-y}{2} \Big\|^p
\end{equation*}
can be shown to be equivalent with $E$ being a UMD-space (see \cite{BURKHOLDEREXP}).
These ideas have been applied (again in \cite{BURKHOLDEREXP}) to prove the boundedness of Doob's maximal operator 
and we will now use them to study the Rademacher maximal function. More precisely, we will
show that for a fixed $p\in (1,\infty )$, a constant $C$ is such that 
$\E | X_R^* |^p \leq C \| X \|_p^p$ for all finite simple
martingales $X=(X_j)_{j=1}^N$ in $\mathcal{X}$ 
if and only if there exists a suitable majorant for the real-valued function
\begin{equation*}
  u(\mathcal{T},T) = \mathcal{R}(\mathcal{T})^p - C \| T \|^p,
\end{equation*}
defined for finite subsets $\mathcal{T}$ of operators in $\mathcal{X}$ and $T\in \mathcal{X}$. Observe that
$\E | X_R^* |^p - C \| X \|_p^p \leq 0$ can equivalently be written as
\begin{equation*}
  \E u \Big( \{ X_j \}_{j=1}^N ,  X_N \Big) \leq 0 ,
\end{equation*}
since $\| X \|_p^p = \E \| X_N \|^p$.

\begin{prop}
\label{concaveprop}
The estimate 
\begin{equation*}
  \E u \Big( \{ X_j \}_{j=1}^N ,  X_N \Big) \leq 0
\end{equation*}
holds for all finite simple martingales $X=(X_j)_{j=1}^N$ in $\mathcal{X}$ 
if and only if there exists a function $v$ satisfying
\begin{enumerate}
  \item $v(\mathcal{T},T) \geq u(\mathcal{T},T)$
  \item $v(\{ T \} , T ) \leq 0$
  \item $v(\mathcal{T} \cup \{ T \} , T ) = v(\mathcal{T},T)$
  \item $v(\mathcal{T}, \cdot )$ is concave
\end{enumerate}
for all finite subsets $\mathcal{T}$ of $\mathcal{X}$ and all $T\in \mathcal{X}$.
\end{prop}

The proof of sufficiency is based on the following lemma.

\begin{lemma}
\label{concavelemma}
  Suppose that $v$ is as in Proposition \ref{concaveprop} 
  and that $(X_j)_{j=1}^N$ is a simple martingale in $\mathcal{X}$.  
  Then, for all $2\leq k \leq N$, we have
  \begin{equation*}
    \E v \Big( \{ X_j \}_{j=1}^k , X_k \Big)  
    \leq \E v \Big( \{ X_j \}_{j=1}^{k-1} , X_{k-1} \Big) .
  \end{equation*}
  \begin{proof}
    Let us fix a $k$ and write $\mathcal{F}_j$ for the $\sigma$-algebra generated by $X_1, \ldots , X_j$.
    By the simplicity of $(X_j)_{j=1}^N$, the set $\{ X_j \}_{j=1}^{k-1}$ has a finite number $s$ of different
    possibilities $\mathcal{T}_1, \ldots , \mathcal{T}_s \subset \mathcal{X}$ so that the 
    event $A_r$ of $\mathcal{T}_r$ happening is in $\mathcal{F}_{k-1}$.
    Now, using the third property of $v$ we get 
    \begin{equation*}
      v \Big( \{ X_j \}_{j=1}^k , X_k \Big)
      = v \Big( \{ X_j \}_{j=1}^{k-1} \cup \{ X_k \} , X_k \Big) 
      = v \Big( \{ X_j \}_{j=1}^{k-1} , X_k \Big) 
      = \sum_{r=1}^s 1_{A_r} v (\mathcal{T}_r , X_k )
    \end{equation*}
and so the fourth property with the aid of Jensen's inequality implies
    \begin{equation*}
      \E \Big( v(\mathcal{T}_r,X_k) \Big| \mathcal{F}_{k-1} \Big) 
      \leq v \Big( \mathcal{T}_r , \E ( X_k | \mathcal{F}_{k-1} ) \Big)
      = v(\mathcal{T}_r, X_{k-1}) .
    \end{equation*}
    Thus
    \begin{align*}
      \E v \Big( \{ X_j \}_{j=1}^k , X_k \Big) 
      &= \sum_{r=1}^s \E \Big( 1_{A_r} v (\mathcal{T}_r, X_k) \Big) \\
      &= \sum_{r=1}^s \E \Big( 1_{A_r} \E \Big( v(\mathcal{T}_r,X_k) \Big| \mathcal{F}_{k-1} \Big) \Big) \\
      &\leq \sum_{r=1}^s \E \Big( 1_{A_r} v(\mathcal{T}_r, X_{k-1}) \Big) \\
      &= \E v \Big( \{ X_j \}_{j=1}^{k-1} , X_{k-1} \Big) ,
    \end{align*}
    where the second equality relies on $A_r$'s belonging to $\mathcal{F}_{k-1}$.
  \end{proof}
\end{lemma}

\begin{proof} (Proof of Proposition \ref{concaveprop}.)

With the aid of the Lemma \ref{concavelemma}, the existence of a desired $v$ is now readily seen to imply that
\begin{equation*}
  \E u \Big( \{ X_j \}_{j=1}^N , X_N \Big) 
  \leq \E v \Big( \{ X_j \}_{j=1}^N , X_N \Big) 
  \leq \E v \Big( \{ X_j \}_{j=1}^{N-1} , X_{N-1} \Big) 
  \leq \ldots \leq \E v \Big( \{ X_1 \} , X_1 \Big) \leq 0.
\end{equation*}

On the other hand, the validity of $\E u \Big( \{ X_j \}_{j=1}^N ,  X_N \Big) \leq 0$ for finite simple martingales
enables us to construct the
auxiliary function $v$ with the desired properties by defining
\begin{equation*}
  v(\mathcal{T},T) = \sup \E u \Big( \{ X_j \}_{j=1}^N \cup \mathcal{T}, X_N \Big) ,
\end{equation*}
where the supremum is taken over all finite and simple martingales $(X_j)_{j=1}^N$ (where $N$ is allowed to vary)
for which $X_1 = T$ almost surely.
Let us check that the required properties are satisfied.
For the first property, take $N=1$ and $X_1 = T$ almost surely to see that
\begin{equation*}
  u(\mathcal{T},T) = \mathcal{R}(\mathcal{T})^p - C \| T \|^p
  \leq \mathcal{R}(\mathcal{T} \cup \{ T \} )^p - C \| T \|^p
  = \E \Big( \mathcal{R}(\mathcal{T} \cup \{ X_1 \} )^p - C \| X_1 \|^p \Big) 
  \leq v(\mathcal{T},T).
\end{equation*}
For the third one, it suffices to note that if $X_1 = T$ almost surely, then $\{ T \} \subset \{ X_j \}_{j=1}^N$
almost surely and so $v(\mathcal{T} \cup \{ T \} , T) = v(\mathcal{T},T)$.
The second property follows from the assumption and the third property:
Let $X=(X_j)_{j=1}^N$
be a simple martingale with $X_1 = T$ almost surely. Now
\begin{equation*}
  \E u \Big( \{ X_j \}_{j=1}^N \cup \emptyset , X_N \Big) \leq 0
\end{equation*}
and so $v(\emptyset , T ) \leq 0$. By the third property,
\begin{equation*}
  v(\{ T \} , T ) = v(\emptyset , T) \leq 0.
\end{equation*}

To see that $v(\mathcal{T}, \cdot )$ is concave, take operators $T_1$ and $T_2$ and put 
$T = \alpha T_1 + (1-\alpha )T_2$ for some $0 < \alpha < 1$. We need to show that
$v(\mathcal{T} , T) \geq \alpha v(\mathcal{T} , T_1) + (1-\alpha ) v(\mathcal{T} , T_2)$. To do this, take
$m_1$ and $m_2$ such that $m_i < v(\mathcal{T} , T_i)$. Now there exist finite simple martingales
$(X_j^{(i)})_{j=1}^N$ (defined on the unit interval) such that $X_1^{(i)} = T_i$ almost surely and
\begin{equation*}
  \E u\Big( \{ X_j^{(i)} \}_{j=1}^N \cup \mathcal{T} , X_N^{(i)} \Big) > m_i . 
\end{equation*}
Let $X_1 = T$ almost surely and define
\begin{equation*}
  X_j(t) =
  \begin{cases}
    X_{j-1}^{(1)}(\frac{t}{\alpha}) , \quad t\in [0,\alpha ) \\
    X_{j-1}^{(2)}(\frac{t - \alpha}{1 - \alpha}) , \quad t\in [\alpha , 1)
  \end{cases}
\end{equation*}
for $j=2,\ldots , N+1$.

\begin{figure}[h!]
\centering
  \setlength{\unitlength}{1bp}%
    \begin{picture}(250.40, 86.04)(0,0)
  \put(0,0){\includegraphics{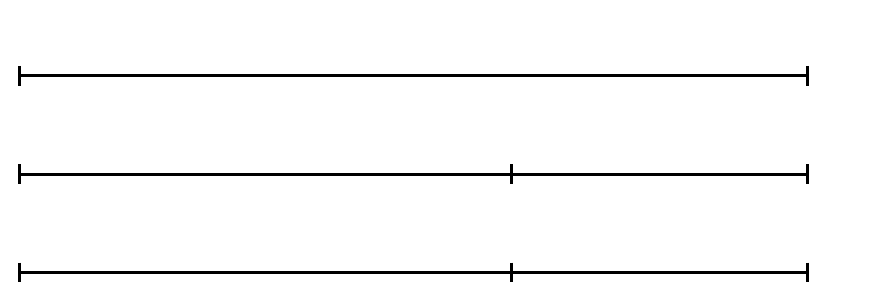}}
  \put(99.21,73.70){\fontsize{8.54}{10.24}\selectfont $X_1 = T$}
  \put(39.69,45.35){\fontsize{8.54}{10.24}\selectfont $X_2 = X_1^{(1)} = T_1$}
  \put(153.07,45.35){\fontsize{8.54}{10.24}\selectfont $X_2 = X_1^{(2)} = T_2$}
  \put(53.86,17.01){\fontsize{8.54}{10.24}\selectfont $X_3 = X_2^{(1)}$}
  \put(164.41,17.01){\fontsize{8.54}{10.24}\selectfont $X_3 = X_2^{(2)}$}
  \end{picture}%
\caption{The construction of $X_1$, $X_2$ and $X_3$}
\end{figure}
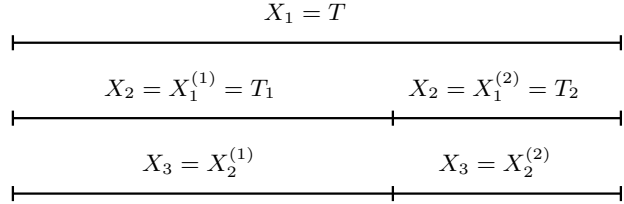

A moments reflection assures us that $(X_j)_{j=1}^{N+1}$ is also a simple martingale. Now
\begin{align*}
  v(\mathcal{T},T) &> \E u \Big( \{ X_j \}_{j=1}^{N+1} \cup \mathcal{T} , X_{N+1} \Big) \\
  &\geq \E u \Big( \{ X_j \}_{j=2}^{N+1} \cup \mathcal{T} , X_{N+1} \Big) \\
  &= \int_0^{\alpha} 
  u \Big( \Big\{ X_j^{(1)}(\frac{t}{\alpha}) \Big\}_{j=1}^N \cup \mathcal{T} , 
  X_N^{(1)}(\frac{t}{\alpha}) \Big) \D t \\
  &+ \int_{\alpha}^1
  u \Big( \Big\{ X_j^{(2)}(\frac{t-\alpha}{1-\alpha}) \Big\}_{j=1}^N \cup
  \mathcal{T} , X_N^{(2)}(\frac{t-\alpha}{1-\alpha}) \Big) \D t \\
  &= \alpha \int_0^1 
  u \Big( \{ X_j^{(1)}(s) \}_{j=1}^N \cup \mathcal{T} , X_N^{(1)}(s) \Big) \D s
  + (1-\alpha ) \int_0^1
  u \Big( \{ X_j^{(2)}(s) \}_{j=1}^N \cup 
  \mathcal{T} , X_N^{(2)}(s) \Big) \D s \\
  &> \alpha m_1 + (1-\alpha ) m_2 .
\end{align*}
Letting $m_i \to v(\mathcal{T} , T_i)$ we get concavity. The proof of Proposition \ref{concaveprop} is now complete.
\end{proof}

\begin{remark}
Had we assumed in Proposition \ref{concaveprop} that
$\E u \Big( \{ X_j \}_{j=1}^N ,  X_N \Big) \leq 0$ holds only for standard Haar martingales, 
we would have obtained a function $v$ satisfying properties (1)-(3) but for which
$v(\mathcal{T},\cdot )$ is only midpoint concave. 
Indeed, suppose that the supremum in the
definition of $v$ is taken over finite standard Haar martingales and observe that properties other than concavity
follow exactly as above. In the proof of midpoint concavity, let
$T = ( T_1 + T_2 ) / 2$ and define $(X_j)_{j=1}^{2N+1}$ as follows:
\begin{align*}
  X_1 &= T \quad \text{almost surely}, \\
  X_2(t) &=
  \begin{cases}
    X_1^{(1)}(2t) = T_1 , \quad t\in [0,1/2 ) , \\
    X_1^{(2)}(2t-1) = T_2 , \quad t\in [1/2 , 1) ,
  \end{cases} \\
  X_{2j-1} (t) &=
  \begin{cases}
    X_j^{(1)}(2t) , \quad t\in [0,1/2 ) , \\
    X_{2j-2}(t) , \quad t\in [1/2 , 1) ,
  \end{cases} \\
  X_{2j} (t) &=
  \begin{cases}
    X_{2j-1}(t) , \quad t\in [0,1/2 ) , \\
    X_j^{(2)}(2t-1) , \quad t\in [1/2 , 1) .
  \end{cases}
\end{align*}

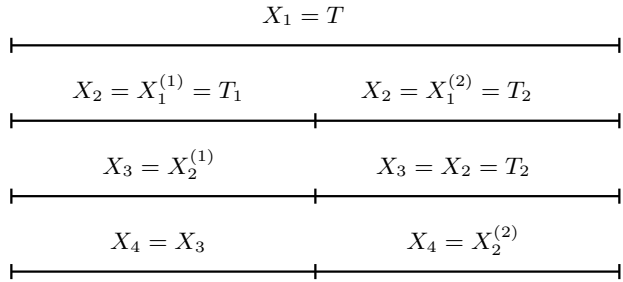
\begin{figure}[h!]
\centering
  \setlength{\unitlength}{1bp}%
  \begin{picture}(238.11, 114.38)(0,0)
  \put(0,0){\includegraphics{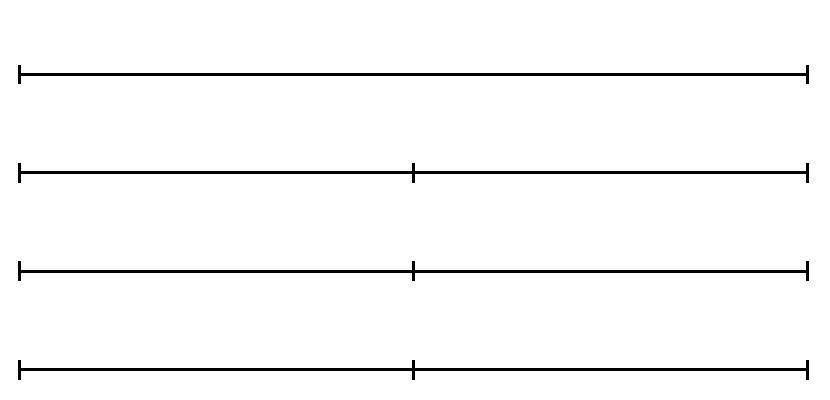}}
  \put(99.21,102.05){\fontsize{8.54}{10.24}\selectfont $X_1 = T$}
  \put(136.06,73.70){\fontsize{8.54}{10.24}\selectfont $X_2 = X_1^{(2)} = T_2$}
  \put(39.69,45.35){\fontsize{8.54}{10.24}\selectfont $X_3 = X_2^{(1)}$}
  \put(141.73,45.35){\fontsize{8.54}{10.24}\selectfont $X_3 = X_2 = T_2$}
  \put(42.52,17.01){\fontsize{8.54}{10.24}\selectfont $X_4 = X_3$}
  \put(153.07,17.01){\fontsize{8.54}{10.24}\selectfont $X_4 = X_2^{(2)}$}
  \put(28.35,73.70){\fontsize{8.54}{10.24}\selectfont $X_2 = X_1^{(1)} = T_1$}
  \end{picture}%
\caption{The construction of $X_1$, $X_2$, $X_3$ and $X_4$}
\end{figure}

This way $(X_j)_{j=1}^{2N+1}$ becomes a standard Haar martingale and calculations similar as in the proof of
Proposition \ref{concaveprop} give us
$v(\mathcal{T} , T) \geq v(\mathcal{T} , T_1) / 2 + v(\mathcal{T} , T_2) / 2$. 
\end{remark}

In conclusion, we state:

\begin{theorem}
\label{concave}
  Let $1 < p < \infty$. Then
  $\mathcal{X}$ has $\text{RMF}_p$ if and only if there exists a 
  function $v$ such that for some constant $C$,
  \begin{enumerate}
    \item $v(\mathcal{T},T) \geq \mathcal{R}(\mathcal{T})^p - C \| T \|^p$,
    \item $v(\{ T \} , T ) \leq 0$,
    \item $v(\mathcal{T} \cup \{ T \} , T ) = v(\mathcal{T},T)$,
    \item $v(\mathcal{T}, \cdot )$ is midpoint concave,
  \end{enumerate}
  for all finite subsets $\mathcal{T}$ of $\mathcal{X}$ and all $T\in \mathcal{X}$.
  \begin{proof}
    If $\mathcal{X}$ has $\text{RMF}_p$, there exists a constant $C$ is such that 
    $\E | X_R^* |^p \leq C \| X \|_p^p$ especially for all standard Haar martingales $X=(X_j)_{j=1}^N$ in $\mathcal{X}$.
    Equivalently, 
    \begin{equation*}
      \E \Big( \mathcal{R} \Big( X_j : 1 \leq j \leq N \Big)^p - C \| X_N \|^p \Big) \leq 0
    \end{equation*}
    for standard Haar martingales $X=(X_j)_{j=1}^N$, which by Proposition
    \ref{concaveprop} enables us to construct a desired $v$.

    To show the converse, we first sketch a proof of the known fact that
    midpoint concave functions that are locally bounded from below are actually
    concave. 
    Suppose that a function $f:\mathcal{X} \to \R$ is midpoint concave but not concave. 
    Then there exist $T_0,T_1\in\mathcal{X}$ such that even though
    \begin{equation*}
      f((1-\alpha )T_0 + \alpha T_1) \geq (1-\alpha ) f(T_0) + \alpha f(T_1)
    \end{equation*}    
    holds (by induction from midpoint concavity) for all $\alpha$ of the dyadic form
    $m2^{-k}$ with $k\geq 1$ and $m\in \{ 1,\ldots , 2^k \}$, it does not hold for some 
    $\alpha' \in (0,1)$.
    Assuming with no loss of generality that $f(T_1) \geq f(T_0)$, 
    we claim that such an $f$ can not be locally 
    bounded from below. Let us write $T_{\alpha} = (1 - \alpha ) T_0 + \alpha T_1$
    and $c_{\alpha} = (1- \alpha) f(T_0) + \alpha f(T_1)$ for $\alpha \in (0,1)$ so that 
    $\delta = c_{\alpha'} - f(T_{\alpha'}) > 0$. 
    One can now express $\alpha'$ as the midpoint
    of an interval $[\alpha , \alpha''] \subset (0,1)$,
    where $\alpha$ is dyadic and so close to $\alpha'$ that
    \begin{equation*}
      f(T_{\alpha}) - f(T_{\alpha'}) \geq (1 - \alpha ) f(T_0) + \alpha f(T_1) - 
      f(T_{\alpha'}) \geq \delta / 2 .
    \end{equation*} 
    By midpoint concavity of $f$ we have 
    $f(T_{\alpha'})\geq (f(T_{\alpha}) + f(T_{\alpha''})) / 2$ and so
    \begin{equation*}
      f(T_{\alpha''}) \leq f(T_{\alpha'}) - (f(T_{\alpha}) - f(T_{\alpha'})) 
      \leq c_{\alpha'} - 3\delta / 2 \leq c_{\alpha''} - 3\delta / 2 ,
    \end{equation*}    
    where the last inequality follows from $\alpha'' \geq \alpha'$ by our
    assumption $f(T_1) \geq f(T_0)$.
    Hence, starting with an $\alpha' \in (0,1)$ such that
    $f(T_{\alpha'}) \leq c_{\alpha'} - \delta$ we find an $\alpha'' \in (0,1)$ for which
    $f(T_{\alpha''}) \leq c_{\alpha''} - 3\delta / 2$.
    Continuing this way we see that $f$ can not be locally bounded from below since the numbers
    $c(\alpha)$ are bounded from above by $f(T_1)$.

    Suppose then that there exists a function $v$ with the listed properties.
    By the first property, $v(\mathcal{T} , \cdot )$ is locally bounded from below. 
    Hence $v$ is concave and by Proposition \ref{concaveprop} we have
    $\E | X_R^* |^p \leq C \| X \|_p^p$ especially for all finite simple martingales $X=(X_j)_{j=1}^N$.
    By Theorem \ref{rmffilt} (or just by Lemma \ref{finitealgs}) $\mathcal{X}$ has $\text{RMF}_p$.
  \end{proof}
\end{theorem}

Observe that this is another way to see that to have
the condition $\E | X_R^* |^p \leq C \| X \|_p^p$ for finite simple martingales it suffices 
to check it for standard Haar martingales.

\section*{Acknowledgements}

I would like to express my gratitude to The Finnish Centre of Excellence in Analysis and Dynamics Research
and to The Finnish National Graduate School in Mathematical Analysis and Its Applications for their support.
Many thanks to the referee who carefully read the manuscript and rightfully suggested some corrections.
This article is a part of my Licentiate Thesis 
which I have written under the supervision of Tuomas Hytönen at the University of Helsinki.

\bibliographystyle{plain}
\bibliography{viitteet}

\end{document}